\documentclass[11pt,fleqn]{article}
\usepackage[bookmarks=false,pdfpagemode=UseOutlines ,plainpages=false
,hypertexnames=false ,pdfpagelabels ,hyperindex=true,colorlinks=true,citecolor=blue]{hyperref}

\RequirePackage[OT1]{fontenc}
\RequirePackage{amsthm,amsmath}
\RequirePackage[round]{natbib}

\usepackage{amsmath,amsfonts,amssymb,euscript,mathrsfs}
\usepackage{slashbox}
\usepackage{float}
\usepackage[final]{pdfpages}
\usepackage{dsfont}
\usepackage{graphicx}
\usepackage{txfonts}
\usepackage[T1]{fontenc} 
\usepackage{fullpage}

\def\1{\mbox{\bf 1}} 
\def\R{\mathbb{R}}

\def\N{\mathbb{N}}
\def\P{\mathbb{P}}
\def\E{\mathbb{E}}
\def\L{\mathbb{L}}

\def\R{\mathbb{R}}

\def\Z{\mathbb{Z}}

\def\v{\mbox{Var\,}}

\def\cov{\mbox{Cov}}

\def\ep{\varepsilon}
\newcommand{\nor}[1]{\left|#1\right|_{\infty,\epsilon}}

\newtheorem{theo}{Theorem}
\newtheorem{lem}{Lemma}

\newtheorem{cor}{Corollary}

\begin{document}
\author{Lionel Truquet \footnote{UMR 6625 CNRS IRMAR, University of Rennes 1, Campus de Beaulieu, F-35042 Rennes Cedex, France and}
\footnote{Campus de Ker-Lann, rue Blaise Pascal, BP 37203, 35172 Bruz cedex, France. {\it Email: lionel.truquet@ensai.fr}.}
 }
\title{Root$-n$ consistent estimation of the marginal density in some time series models}
\date{}
\maketitle

\begin{abstract}
In this paper, we consider the problem of estimating the marginal density in some nonlinear autoregressive time series models for which the conditional mean and variance have a parametric specification.
Under some regularity conditions, we show that a kernel type estimate based on the residuals can be root$-n$ consistent even if the noise density is unknown. 
Our results, which are shown to be valid for classical time series models such as ARMA or GARCH processes, extend substantially the existing results obtained for some homoscedatic time series models. Asymptotic expansion of our estimator is obtained by combining some martingale type arguments and a coupling method for time series which is of independent interest.
We also study the uniform convergence of our estimator on compact intervals. 
\end{abstract}

\section{Introduction}
Estimating the marginal density of a stationary time series has been extensively studied in the literature.
Kernel density estimation is probably one of the most popular methods used for this problem and the properties of 
the so-called Parzen-Rosenblatt estimator have been investigated under various mixing type conditions. See for instance \citet{Rob}, \citet{Dou1},\citet{Dou2} or \citet{Rou}. See also the monograph of \citet{Bosq} for kernel density estimation for strong mixing sequences and \citet{Ded} for numerous weak dependence conditions ensuring consistency properties of this estimator.  
  
However, when additional structure is assumed for the stochastic process of interest, kernel density estimation can be used more cleverly 
for getting sharper rates of convergence, in particular $\sqrt{n}-$consistency. This atypical rate of convergence in nonparametric density estimation has been first noticed for the estimation of the density of some functionals of independent random variables. See \citet{Frees}, \citet{SW3} and \citet{Gine}.    
In time series, existing contributions exploits the representation of the marginal density as a convolution product between the innovation density and the marginal density of a predictable process. Such an approach has been used by \citet{Cao}, \citet{SW2} and \citet{SW1} for estimating the marginal density of invertible moving average processes. In the last contribution, sharper results are obtained for possibly infinite moving averages processes.  
More recently, \citet{Wu} obtained some results for nonlinear and homoscedastic autoregressive processes of order $1$ for which the conditional mean has a parametric specification. Another recent contribution has been made by \citet{Del} who constructed a $\sqrt{n}-$consistent estimator of the density of the log-volatility for a GARCH$(1,1)$ process. Note however that the purpose of this latter contribution is not the estimation of the marginal distribution and the volatility process is not directly observed. Moreover, the approach used seems specific to the autoregressive equation followed by the GARCH$(1,1)$.  

In the literature, $\sqrt{n}-$consistent estimation of the marginal density in conditionally heteroscedastic time series models has not been considered. Moreover, even in the homoscedastic case, a general approach has not been studied for getting this convergence rate. In this paper, we consider the problem of estimating the marginal density with the $\sqrt{n}$ rate of convergence in some autoregressive time series models, conditionally homoscedastic or heteroscedastic. We will restrict our study to short memory models with a location-scale formulation 
$$X_t=m_t(\theta)+\sigma_t(\theta)\varepsilon_t,\quad t\in\Z.$$
where the conditional mean $m_t(\theta)$ and volatility $\sigma_t(\theta)$ depends smoothly on a finite-dimensional parameter $\theta$.
With respect to the existing results, our approach covers lots of cases, from the ARMA processes with independent and identically distributed innovations to ARMA processes with a GARCH noise.    
Let us also mention that our contribution gives an answer to a question addressed in \citet{Zh}, a paper in which an estimator similar to our was suggested for density estimation in autoregressive time series models. 
However, apart from some classical smoothness conditions, the root-$n$ consistency of our estimator is only guaranteed under the square integrability, with respect to the noise distribution, of the conditional density of the marginal $X_t$ given the noise component $\varepsilon_t$. This condition is not always satisfied and has to be checked for the model under study. In this paper, we show in particular that ARMA processes with GARCH errors satisfy in general this integrability condition. 
A similar integrability condition has been observed by \citet{Mul} for estimating the marginal density in some homoscedastic regression models. See also \citet{SWMetrika} who showed that estimating a convolution of some powers of independent random variables can lead to a slower rate of convergence when this conditions fails to hold.      

The paper is organized as follows. In Section $2$, we define our estimator and give its asymptotic properties. 
In Section $3$, we check the assumptions of our Theorems for some standard examples of time series models. We also compare our assumptions with that used in the aforementioned references. Section $4$ is devoted to a comparison by simulation of our estimator with the standard
Parzen-Rosenblatt estimator. Proofs of our results are postponed to the last section of the paper.     

\section{Marginal density estimation of a time series}
\subsection{Model and marginal density estimator}
We first introduce the general model used in the sequel.
Let $\left(\ep_t\right)_{t\in\Z}$ be a sequence of i.i.d square integrable random variables.
If $\Theta$ denotes a Borel subset of $\R^d$, we consider two measurable functions $H,G:\Theta\times \R^{\N}\rightarrow \R$.
We assume that for a $\theta_0\in\Theta$, $(X_t)_{t\in\Z}$ is a stationary process such that
\begin{equation}\label{crucial}
X_t=H\left(\theta_0;X_{t-1},X_{t-2},\ldots\right)+\ep_tG(\theta_0;X_{t-1},X_{t-2},\ldots).
\end{equation}
Note that the two functions $H$ and $G$ will be more precisely defined $\lambda_d\otimes \P_X$ almost everywhere, where 
$\lambda_d$ denotes the Lebesgue measure on $\R^d$ and $\P_X$ the probability distribution of $(X_{t-1},X_{t-2},\ldots)$. 
We also assume that $X_t\in\sigma\left(\ep_t,\ep_{t-1},\ldots\right)$, i.e
$$X_t=E\left(\ep_t,\ep_{t-1},\ldots\right),$$
for a suitable measurable function $E:\R^{\N}\rightarrow \R$ defined $\P_{\varepsilon}$ almost everywhere. 
We also set 
$$m_t(\theta)=H\left(\theta;X_{t-1},X_{t-2},\ldots\right),\quad \sigma_t(\theta)=G\left(\theta;X_{t-1},X_{t-2},\ldots\right).$$
The realizations of all the past values are not available. Then we assume that there exist measurable functions $H_t,G_t:\Theta\times \R^t\rightarrow \R$ such that $H_t\left(\theta;X_{t-1},\ldots,X_1\right)$ (resp. $G_t\left(\theta;X_{t-1},\ldots,X_1\right)$) is an approximation of $m_t(\theta)$ (resp. $\sigma_t(\theta)$). Then we use the notations
$$\overline{m}_t(\theta)=H_t\left(\theta;X_{t-1},\ldots,X_1\right),\quad \overline{\sigma}_t(\theta)=G_t\left(\theta;X_{t-1},\ldots,X_1\right).$$
Our estimator is based on the representation of the marginal density $f_X$ of the stationary process $(X_t)_{t\in\Z}$ as a smooth functional of the noise density $f_{\ep}$. 
More precisely, setting $X_t^{-}=\left(X_{t-1},X_{t-2},\ldots\right)$ and denoting by $f\left(\cdot\vert X_t^{-}\right)$ denotes the conditional density of $X_t$ given $X_t^{-}$, we have for $v\in\R$,
\begin{equation}\label{departure}
$$f_X(v)=\E\left[f\left(v\vert X_t^{-}\right)\right]=\E\left[\frac{1}{\sigma_i(\theta_0)}f_{\ep}\left(\frac{v-m_i(\theta_0)}{\sigma_i(\theta_0)}\right)\right].
\end{equation}
Imagine first that a sample $\left(X_i,m_i(\theta_0),\sigma_i(\theta_0)\right)_{1\leq i\leq n}$ is available.
Then the vector of innovations $\left(\ep_1,\ldots,\ep_n\right)$ is also observed. The noise density $f_{\ep}$ can be estimated by the classical Parzen-Rosenblatt kernel estimator. If $K:\R\rightarrow\R_+$ be a probability density with compact support $[-1,1]$, which will be assumed to be continuously differentiable in the sequel, we set 
$$\hat{f}_{\ep}(z)=\frac{1}{n}\sum_{i=1}^n K_b\left(z-\ep_i\right),\quad K_b(x)=\frac{1}{b}K\left(\frac{x}{b}\right).$$
Then, using the expression (\ref{departure}), we define the following unfeasible estimator
\begin{eqnarray*}
\check{f}_X(v)&=&\frac{1}{n}\sum_{i=1}^n \frac{1}{\sigma_i(\theta_0)}\hat{f}_{\ep}\left(L_{v,i}(\theta_0)\right)\\ 
&=&\frac{1}{n^2}\sum_{i,j=1}^n\frac{1}{\sigma_i(\theta_0)}K_b\left[L_{v,i}(\theta_0)-\ep_j\right],
\end{eqnarray*}
with $L_{v,i}(\theta)=\frac{v-m_i(\theta)}{\sigma_i(\theta)}$ for $(v,\theta)\in \R\times\Theta$. 
In practice, the parameter $\theta_0$ has to be estimated and only the vector $(X_1,X_2,\ldots,X_n)$ is observed. 
Let us introduce additional notations. For $(v,\theta)\in \R\times\Theta$, we set $\ep_j(\theta)=\frac{X_j-m_j(\theta)}{\sigma_j(\theta)}$. We call the process $\left(\ep_j(\cdot)\right)_{j\in\Z}$ the residual process.  
We also denote by $\overline{\ep}_j(\theta)$ and $\overline{L}_{v,i}(\theta)$ the truncated versions of $\ep_j(\theta)$ and $L_{v,i}(\theta)$ respectively, e.g $\overline{\ep}_j(\theta)=\frac{X_j-\overline{m}_j(\theta)}{\overline{\sigma}_j(\theta)}$.

Then, if $\hat{\theta}$ denotes an estimator of $\theta_0$, the feasible estimator of $f_X(v)$ is defined by
$$\hat{f}_X(v)=\frac{1}{n^2}\sum_{i,j=1}^n\frac{1}{\overline{\sigma}_i(\hat{\theta})}K_b\left[\overline{L}_{v,i}(\hat{\theta})-\overline{\ep}_j(\hat{\theta})\right].$$
Note that $\left(\overline{\ep}_j(\hat{\theta})\right)_{1\leq j\leq n}$ are the residuals obtained after the estimation step.

In the homosecedastic case, i.e there exists $\sigma>0$ such that $\sigma_t(\theta)=\sigma$ for all $(t,\theta)\in\Z\times \Theta$,
our estimator is simply defined by 
\begin{equation}\label{reduction}
\hat{f}_X(v)= \frac{1}{n^2}\sum_{i,j=1}^nK_b\left[v-\overline{m}_i(\hat{\theta})-X_j+\overline{m}_j(\hat{\theta})\right].
\end{equation}
Note that the estimation of the variance $\sigma^2$ is unnecessary in the homoscedastic case. Estimator of type (\ref{reduction}) already appears in the literature but using a convolution approach. See for instance \citet{SW2} for linear processes, \citet{Mul} for homoscedastic regression models and \citet{Wu} for some non linear conditionally homoscedastic time series. In this case, the kernel $K$ is a convolution product of type $k*k$ and the estimator (\ref{reduction}) is obtained as a convolution product of two kernel estimators: the Parzen-Rosenblatt estimator, with kernel $k$, of the density of $m_t(\theta_0)$ and that of $f_{\ep}$ with the same kernel. In this paper, we will consider an arbitrary continuously differentiable kernel $K$ and the homoscedastic case as a special case of the conditionally heteroscedastic case, by setting in this case the two quantities  $\sigma_t$ and $\overline{\sigma}_t$ to $1$ in all our statements.     

\subsection{Assumptions and asymptotic behavior of the marginal density estimate}
We now give our assumptions for deriving the asymptotic behavior of the unfeasible estimator $\check{f}_X$ and the feasible estimator $\hat{f}_X$. In the sequel, we will denote by $\Vert\cdot\Vert$ a norm on $\R^d$ whatever the value of the integer $d$. 
We will still denote by $\Vert\cdot\Vert$ the corresponding operator norm. 
For a family $\left\{A(\theta);\theta\in\Theta\right\}$ of matrices and a family $\left\{B(\theta): \theta\in\Theta\right\}$ of real numbers, we set $\nor{A}=\sup_{\theta\in\Theta_{0,\epsilon}}\Vert A(\theta)\Vert$ and $\nor{A,B}=\sup_{\theta,\theta'\in \Theta_{0,\epsilon}}\Vert \frac{A(\theta)}{B(\theta')}\Vert$, where $\Theta_{0,\epsilon}=\left\{\theta\in\Theta: \Vert\theta-\theta_0\Vert<\epsilon\right\}$.
  
Finally, since for $i\in\Z$, $m_i(\theta)$ and $\sigma_i(\theta)$ are measurable functions of $Y_i=\left(\ep_i,\ep_{i-1},\ldots\right)$, we define some coupling versions of these two quantities. For an integer $\ell\geq 1$, we denote by $\left(\ep_j^{(\ell)}\right)_{j\in\Z}$ a copy of $\left(\ep_j\right)_{j\in\Z}$ and we denote by $m_{i\ell}(\theta)$ and $\sigma_{i\ell}(\theta)$ the two random variables 
defined as $m_i(\theta)$ and  $\sigma_i(\theta)$ but for which $Y_i$ is replaced with 
$$Y_{i\ell}=\left(\ep_i,\ep_{i-1},\ldots,\ep_{i-\ell+1},\ep^{(\ell)}_{i-\ell},\ep^{(\ell)}_{i-\ell-1},\ldots\right).$$
One can note that $\left(m_i(\theta),\sigma_i(\theta)\right)$ has the same distribution than $\left(m_i(\theta),\sigma_i(\theta)\right)$.
The interest of such coupling method will be explained in Section $5$.
The following assumptions will be needed.
\begin{description}
\item[A1] The parameter $\theta\in\Theta$ where $\Theta$ is a compact set of $\R^d$. 
\item[A2]
The volatility is bounded away from zero, i.e there exists $\gamma>0$ such that $\inf_{\theta\in\Theta}\sigma_i(\theta)\geq \gamma$ a.s.
We also assume $\inf_{\theta\in\Theta}\overline{\sigma}_i(\theta)\geq \gamma$ a.s.
Moreover, there exists $s,a\in(0,1)$ and $\kappa>0$ such that
$$\E\left[\sup_{\theta\in\Theta_{0,\epsilon}}\left|m_t(\theta)\right|^s+\sup_{\theta\in\Theta_{0,\epsilon}}\left|\sigma_t(\theta)\right|^s\right]<\infty,$$
$$\E\left[\sup_{\theta\in\Theta_{0,\epsilon}}\left|m_i(\theta)-m_{i\ell}(\theta)\right|^s+\sup_{\theta\in\Theta_{0,\epsilon}}\left|\sigma^2_i(\theta)-\sigma^2_{i\ell}(\theta)\right|^s\right]\leq \kappa a^{\ell}.$$
and 
$$\E\left[\sup_{\theta\in\Theta_{0,\epsilon}}\left|m_i(\theta)-\overline{m}_i(\theta)\right|^s+\sup_{\theta\in\Theta_{0,\epsilon}}\left|\sigma^2_i(\theta)-\overline{\sigma}^2_i(\theta)\right|^s\right]\leq \kappa a^{i}.$$
 
\item[A3]
The two applications $\theta\mapsto \sigma_t(\theta)$ and $\theta\mapsto m_t(\theta)$ are a.s two times differentiable over $\Theta$. Moreover, there exists $\epsilon>0$ such that the following random variables are integrable.
$$\nor{\dot{m}_i,\sigma_i}^3,\quad \nor{\dot{\sigma_i^2},\sigma^2_i}^3,\quad \nor{\dot{m}_i,\sigma_i}^2\cdot\nor{\dot{\sigma_i^2},\sigma^2_i}^2,\quad \nor{\sigma_i,\sigma_i}\cdot\nor{\ddot{m_i},\sigma_i},$$
$$\nor{m_i,\sigma_i}\cdot\nor{\ddot{\sigma_i^2},\sigma^2_i},\quad \nor{\ddot{\sigma_i^2},\sigma_i^2}^{6/5},$$
where for a function $g:\Theta\rightarrow \R$, $\dot{g}$ and $\ddot{g}$ denote the two first derivatives of $g$. 

\item[A4] There exists an estimator $\hat{\theta}$ of $\theta_0$ such that $\hat{\theta}-\theta_0=O_{\P}\left(\frac{1}{\sqrt{n}}\right)$. 

\item[A5]
The noise density $f_{\ep}$ is bounded and the two first derivatives $f'_{\ep}$, $f''_{\ep}$ are bounded.

\item[A6] Let $I$ a compact interval of the real line. For all $v\in I$, we assume that the ratio $\frac{v-m_t(\theta_0)}{\sigma_t(\theta_0)}$ has a density denoted by $h_v$ and we set for $x\in \R$, 
$$g_v(x)=\E\left[\frac{1}{\sigma_t(\theta_0)}\left|\right.\frac{v-m_t(\theta_0)}{\sigma_t(\theta_0)}=x\right]\cdot h_v(x).$$
We assume that the application $(x,v)\mapsto g_v(x)$ is jointly measurable and also that there exists $s_0>0$ such that for all $v\in I$, 
$$\int g_v(x)^2 d\mu(x),\quad d\mu(x)=\sup_{\vert s\vert\leq s_0}f_{\ep}(x+s)dx.$$ 

\item[A7] The envelope function $G$ defined by $G(x)=\sup_{v\in I}g_v(x)$ satisfies $\int G(x)^{2+o}d\mu(x)<\infty$ for some $o\in (0,1)$. Moreover, there exist some constants $\eta,C>0$ such that 
$$N_{[]}\left(\epsilon,\mathcal{G},\L^2(\mu)\right)\leq C \epsilon^{-\eta},$$
where $N_{[]}\left(\epsilon,\mathcal{G},\L^2(\mu)\right)$ denotes the bracketing numbers of the family $\mathcal{G}_I=\left\{g_v: v\in I\right\}$.
\end{description}

\paragraph{Notes}
\begin{enumerate}
\item
Different constants $a,s$ and $\kappa$ can be found for the three bounds given assumptions {\bf A2}. However, we can alway take the minimal value of the exponents $s$, the maximal value of the constant $a$ and the maximal value of the constants $\kappa\geq 1$. So there is no loss of generality in assuming the same constants for the three bounds. 
\item
Assumption ${\bf A2}$ imposes a restriction on the dependence structure of the time series models. These conditions, which are usually  referred as short-memory properties, are satisfied for the standard ARMA or GARCH processes. Roughly speaking, this weak dependence condition means that a perturbation of initial conditions in the data generating process is forgotten exponentially fast. This type of dependence condition is also used by \citet{Zh} or \citet{Wu}.     
\end{enumerate} 

\paragraph{Discussion of the assumption {\bf A6}.} Our results require some regularity conditions for the family of functions $\left\{g_v: v\in I\right\}$. The integrability condition assumed on $g_v$ is necessary for root $n$ consistency, as shown in Theorem \ref{unfeas} stated below. One can also rely this function to the conditional density of $X_t\vert \varepsilon_t$. Indeed, if $h:\R^2\rightarrow \R_+$ is a measurable function, we have, setting for simplicity of notations $m_t=m_t(\theta_0)$ and $\sigma_t=\sigma_t(\theta_0)$, 
\begin{eqnarray*}
\E\left[h(X_t,\varepsilon_t)\right]&=& \E\int h(m_t+\sigma_t x,x)f_{\varepsilon}(x)dx\\
&=& \E\int h\left(v,\frac{v-m_t}{\sigma_t}\right)\frac{1}{\sigma_t}f_{\varepsilon}\left(\frac{v-m_t}{\sigma_t}\right)dv\\
&=& \int\int h(v,x)g_v(x)f_{\varepsilon}(x)dxdv.
\end{eqnarray*}
This shows that $v\mapsto g_v(x)$ can be seen as a version of the conditional density of $X_t$ given that $\varepsilon_t=x$.
In what follows, we discuss alternative expressions for the function $g_v$ as well as some sufficient conditions ensuring the square integrability of $g_v$ required in {\bf A3}. 
\begin{enumerate}
\item
For the homoscedastic model, the volatility $\sigma_t(\cdot)$ equals to a constant $\sigma$. Then if $m_t(\theta_0)$ has a density denoted by $f_m$, we have $g_v(x)=f_m(v-\sigma x)$. In this case, estimation of parameter $\sigma$ will be unnecessary.
\item
For a pure heteroscedastic model, i.e the conditional mean $m_t(\cdot)$ reduced to a constant $m$, assumption {\bf A3} can hold only if $v\neq m$. For $v\neq m$, one can show that for $x\neq 0$,  $g_v(x)=\frac{1}{\vert x\vert}f_{\sigma}\left(\frac{v-m}{x}\right)$.
In this case, we have for $v\neq m$,  
$$\int g_v(x)^2\mu(dx)\leq \frac{\Vert f_{\varepsilon}\Vert_{\infty}}{\vert v-m\vert}\int_0^{\infty}f^2_{\sigma}(y)dy.$$ 
Then if for instance $f_{\sigma^2}(y)=\frac{1}{2\sqrt{y}}f_{\sigma}(\sqrt{y})$ is bounded, the integrability condition 
given in {\bf A3} is guaranteed as soon as $\E\sigma_t(\theta_0)<\infty$, which is not a strong restriction. 
\item 
In the location-scale case with a non degenerate conditional mean, we assume that the distribution of the couple $\left(m_t(\theta_0),\sigma_t(\theta_0)\right)$ has a density denoted by $f_{m,\sigma}$. Then if $v\in\R$, the distribution of the couple $\left(\frac{v-m_t(\theta_0)}{\sigma_t(\theta_0)},\frac{1}{\sigma_t(\theta_0)}\right)$ has a density $\omega$ given by 
$$\omega(x,y)=\frac{1}{y^3}f_{m,\sigma}\left(v-\frac{x}{y},\frac{1}{y}\right).$$
We deduce that 
\begin{equation}\label{express2}
g_v(x)=\int_0^{1/\gamma}\frac{1}{y^2}f_{m,\sigma}\left(v-\frac{x}{y},\frac{1}{y}\right)dy.
\end{equation}
Using Jensen inequality and a change of variables, one can show that 
\begin{equation}\label{maborne}
\int g_v(x)^2 dx\leq \frac{1}{\gamma}\cdot\E\left[\sigma_t(\theta_0)f_{m,\sigma}\left(m_t(\theta_0),\sigma_t(\theta_0)\right)\right].
\end{equation}
Then the integrability condition given in {\bf A3} follows if $f_{\ep}$ is bounded and if\\ $\E\left[\sigma_t(\theta_0)f_{m,\sigma}\left(m_t(\theta_0),\sigma_t(\theta_0)\right)\right]$ is finite.
\end{enumerate}

\paragraph{Discussion of the assumption {\bf A3}.}
Assumption {\bf A3} imposes some moment restrictions and smoothness conditions. In the pure heteroscedastic case, we only have to check integrability of
$$\nor{\dot{\sigma^2_i},\sigma_i^2}^3\quad \nor{\ddot{\sigma_i^2},\sigma_i^2}^{6/5}.$$
In the homoscedastic case, these conditions reduce to the integrability of
$$\nor{m_i}^2,\quad \nor{\dot{m}_i}^3,\quad \nor{\ddot{m_i}}.$$
These moment restrictions are explained by the technique used for the proof of Theorem $2$ given below and which consists in studying the derivative of $\hat{f}$ with respect to $\hat{\theta}$, the estimator of $\theta_0$. For a general heteroscedastic time series model, other conditions could be possible, the   
single requirement is to get the conclusions of the lemmas \ref{propres}, \ref{propres3}, \ref{propres2}.

We now give the asymptotic behavior of our estimates. We first start with the unfeasible estimator $\check{f}_X$. 

\begin{theo}\label{unfeas}
Assume that there exists $\delta\in (0,1)$ such $nb^{2+\delta}\rightarrow\infty$ and $nb^4\rightarrow 0$. 
\begin{enumerate}
\item
Assume that assumptions {\bf A2}, {\bf A5} and {\bf A6} hold true. Then for all $v\in I$, we have
\begin{equation}\label{approxbien}
\sqrt{n}\left[\check{f}_X(v)-f_X(v)\right]=\frac{1}{\sqrt{n}}\sum_{i=1}^n \left[g_v(\ep_i)+\frac{1}{\sigma_i(\theta_0)}f_{\ep}\left(\frac{v-m_i(\theta_0)}{\sigma_i(\theta_0)}\right)-2f_X(v)\right]+o_{\P}(1).
\end{equation}
In particular, for all $v\in I$, we have $\sqrt{n}\left[\check{f}_X(v)-f_X(v)\right]=O_{\P}(1)$ 
\item
If in addition assumption {\bf A7} holds true, the approximation (\ref{approxbien}) is uniform in $v\in I$. In particular 
$$\sqrt{n}\sup_{v\in I}\left\vert \check{f}_X(v)-f_X(v)\right\vert=O_{\P}(1).$$ 
\end{enumerate}
\end{theo}
\paragraph{Note.} The proof of Theorem $1$ relies on the decomposition of a $U-$statistic for which the degenerate part is shown to be negligible under our bandwidth conditions. The bandwidth condition $\sqrt{n}b^2\rightarrow 0$ is a bias condition, the bias of our estimator has to decrease faster than the rate $1/\sqrt{n}$. However, our estimator will be first approximated by a U-statistic involving $\ell-$dependent random variables in order to facilitate the study of its asymptotic behavior.   

In the next result, we compare the asymptotic behavior of the feasible estimator $\hat{f}_X$ with that of the unfeasible one.
 For $\theta\in \Theta$, we denote by $f_{\theta}$ the density of 
$\ep_i(\theta)$. We also set
$$h_{\theta}(v)=\E\left(\frac{1}{\sigma_i(\theta)}f_{\theta}\left(\frac{v-m_i(\theta)}{\sigma_i(\theta)}\right)\right).$$
In the sequel $\dot{h}_{\theta}(v)$ will denote the partial derivative with respect to $\theta$ of the function $(\theta,v)\mapsto h_{\theta}(v)$. By convention, we represent $\dot{h}_{\theta}(v)$ by a column vector. Moreover, $A^T$ will denote the transpose of a matrix $A$. 

\begin{theo}\label{feas}
Assume that there exists $\delta\in (0,1)$ such that $nb^{3+\delta}\rightarrow\infty$ and $nb^4\rightarrow 0$ and that assumptions {\bf A1-A6} hold true.
Then we have 
$$\sqrt{n}\sup_{v\in I}\left\vert \hat{f}_X(v)-\check{f}_X(v)-\dot{h}_{\theta_0}(v)^T\left(\hat{\theta}-\theta_0\right)\right\vert=o_{\P}(1).$$
\end{theo} 
\paragraph{Note.} For the comparison of the two estimators $\hat{f}$ and $\check{f}$, we do not use $U-$statistics 
arguments. We wanted to avoid additional regularity conditions on the function $g_v$ and its local approximation when $\theta\rightarrow \theta_0$, these conditions being difficult to check for practical examples. For the proof of Theorem $2$, we first show that the effect of truncations of $m_i,\sigma_i$ is negligible. Then, we use $\ell-$dependent approximations of these quantities. Finally, we study a Taylor expansion of order $1$ and control the local approximation of the derivatives using martingale tools and integration with respect to the residual process $\theta\mapsto \frac{X_j-m_j(\theta)}{\sigma_j(\theta)}$. 
Then regularity conditions concern exclusively the densities $\theta\mapsto f_{\theta}$ of this residual process.
One can note that the range of bandwidths allowed in Theorem $2$
is reduced with respect to Theorem $1$.

In the next result, we provide the asymptotic distribution of the feasible estimator $\hat{f}_X$ of $f_X$. 
To this end, it is necessary to use a particular representation of the estimator $\hat{\theta}$. Similar representations are used in 
\citet{Zh} or \citet{Wu}.

\begin{description}
\item[A8] There exists a square integrable process $Z_i(\theta_0)=H_{\theta_0}\left(\ep_{i-1},\ep_{i-2},\ldots\right)$ taking values in $\mathcal{M}_{d,\bar{d}}(\R)$, the space of real matrices of size $d\times\bar{d}$ and a measurable function $F:\R\rightarrow \R^{\bar d}$
such that $\E F(\ep_0)=0$, $Y_i(\theta_0)$ and $F(\ep_0)$ are square integrable, $\E \Vert Z_i(\theta_0)-Z_{i\ell}(\theta_0)\Vert^2 \leq \kappa a^{\ell}$ (where $\kappa>0$ and $a\in (0,1)$) and 
$$\sqrt{n}\left(\hat{\theta}-\theta_0\right)=\frac{1}{\sqrt{n}}\sum_{i=1}^n Z_i(\theta_0)F(\ep_i)+o_{\P}(1).$$
\end{description}

For stating our next result, we define for $v\in I$ and $i\in\Z$, 
$$M_{i,v}=g_v(\ep_i)+\frac{1}{\sigma_i(\theta_0)}f_{\ep}\left(\frac{v-m_i(\theta_0)}{\sigma_i(\theta_0)}\right)+\dot{h}_{\theta_0}(v)^T Z_i(\theta_0)F(\ep_i).$$
\begin{cor}\label{loias}
Under the assumptions {\bf A1-A6} and {\bf A8}, the process $\left(\sqrt{n}\left[\hat{f}_X(v)-f_X(v)\right]\right)_{v\in I}$ converges in the sense of finite dimensional distributions towards a centered Gaussian process $\left(W_v\right)_{v\in I}$ such that for $v_1,v_2\in I$,
$$\cov\left(W_{v_1},W_{v_2}\right)=\cov\left(M_{0,v_1},M_{0,v_2}\right)+\sum_{i\geq 1}\left[\cov\left(M_{0,v_1},M_{i,v_2}\right)+\cov\left(M_{0,v_2},M_{i,v_1}\right)\right].$$
Moreover if assumption {\bf A7} also holds, the convergence occurs in $\ell^{\infty}(I)$.
\end{cor}

\paragraph{Notes} 
\begin{enumerate}
\item
In the homoscedastic case, one can check that the covariance structure of the Gaussian process $\left(W_v\right)_{v\in I}$ is the same 
as in Theorem $2$ of \citet{Wu}. Thus, our result can be seen as an extension of the result obtained in their paper, allowing an arbitrary number of lags in the conditional mean and also conditional heteroscedasticity.
\item
The proofs of Theorem \ref{unfeas}, Theorem \ref{feas} and Corollary \ref{loias} extensively make use of the coupling method discussed at the beginning of Section $5$. This method allows the approximation of some processes by $\ell-$dependent processes which has the same marginal distribution. 
\end{enumerate}
\section{Examples}
In this section, we explain how to check the assumptions ${\bf A1-A8}$ for some classical time series models. 
\subsection{Conditionally homoscedastic times series}
Let us assume that 
$$X_t=m_t(\theta_0)+\varepsilon_t,\quad t\in \Z.$$
We remind that $g_v(x)=f_m(v-x)$ where $f_m$ denotes the density of $m_t$.
In this case, assumption ${\bf A6}$ holds true for instance if $f_m$ is bounded. Moreover, assumption ${\bf A7}$ is satisfied if $f_m$
is Lipschitz and $\int \vert f'_{\varepsilon}(x)\vert dx <\infty$. Indeed, the latter condition entails that the measure $\mu$ of assumption ${\bf A7}$ has a finite mass and in this case, the polynomial decay of the bracketing number is classical. See \citet{VW}, Example $19.7$. 

For homoscedastic and autoregressive time series with one lag, \citet{Wu} obtained a root $n$ consistent estimation of the marginal density by using a representation of the density of $m_t(\theta_0)+\varepsilon_t$ as a convolution product.   
In that paper, similar regularity assumptions are used for the noise distribution. These authors use bandwidth conditions similar to ours (see Theorem $2$ of their paper). Their moment conditions for $\theta\mapsto m_t(\theta_0)$ and its derivative are less restrictive than ours but at the same time more regularity conditions on the density of $m_t(\theta_0)$ have to be checked for their non linear models. See Assumption $4$ and Assumption $6$ of that paper for a precise statement of their regularity conditions. One advantage of our approach is to present a unified approach for homoscedastic and heteroscedastic time series and for which the dynamic can depend on an arbitrary and possibly infinite number of lags.

\paragraph{The case of ARMA processes.} Let us now consider the case of ARMA processes, i.e there exist two integers $p$ and $q$ such that
$$X_t=\eta_0+\sum_{i=1}^p a_{0i} \left(X_{t-i}-\eta_0\right)+\varepsilon_t-\sum_{j=1}^q b_{0j} \varepsilon_{t-j},\quad t\in\Z.$$ 
As usual, we assume that for $\theta=\left(\eta,a_1,\ldots,a_p,b_1,\ldots,b_q\right)\in\Theta$, the roots of the two polynomials $\mathcal{P}(z)=1-\sum_{i=1}^pa_i z^i$ and $\mathcal{Q}(z)=1-\sum_{j=1}^q b_j z^j$ are outside the unit disc. 
Then defining 
$$Z_t(\theta)=X_t-\eta-\sum_{j=1}^pa_j\left(X_{t-j}-\eta\right)+\sum_{j=1}^q b_j Z_{t-j}(\theta)$$
and $\underline{Z}_t(\theta)=\left(Z_t(\theta),Z_{t-1}(\theta),\ldots,Z_{t-q+1}(\theta)\right)^T$, we have
$$\underline{Z}_t(\theta)=A_1(\theta)\underline{Z}_{t-1}(\theta)+B_{1,t}(\theta),$$
where $A_1(\theta)$ denotes the companion matrix associated to $a_1,\ldots, a_p$ and 
$$B_{1,t}(\theta)=\left(X_t-\eta-\sum_{j=1}^pa_j(X_{t-j}-\eta),0,\ldots,0\right)^T.$$
Then we have 
\begin{equation}\label{expan}
\underline{Z}_t(\theta)=\sum_{j=0}^{\infty}A_1(\theta)^j B_{1,t-j}(\theta),\\
m_t(\theta)=\eta+\sum_{j=1}^pa_j\left(X_{t-j}-\eta\right)-\sum_{j=1}^q b_j Z_{t-j}(\theta).
\end{equation}
Moreover $\overline{m}_t(\theta)$ can be defined by setting $X_0,X_{-1},\ldots$ to $0$ in (\ref{expan}).
Now we explain how to check assumptions ${\bf A3-A6}$. Assumption ${\bf A2}$ will be checked directly for ARMA-GARCH processes.
\begin{itemize}
\item
Assumption ${\bf A3}$ holds if $\E\vert\varepsilon_t\vert^3<\infty$. Indeed, using the well-known infinite moving-average representation 
$$X_t=\zeta_0+\ep_t+\sum_{j=1}^{\infty}\zeta_j \varepsilon_{t-j},\quad \sum_{j=1}^{\infty}\vert \zeta_j\vert<\infty,$$
we have $\E \vert X_t\vert^3<\infty$. Then assumption {\bf A3} holds true using (\ref{expan}) (the order of the derivative can be arbitrary in this example).
\item
Now assumptions ${\bf A6-A7}$ follows from the fact that $f_m$ is bounded and Lipschitz. Indeed, using the infinite moving average representation, we have $m_t(\theta_0)=\zeta_0+\sum_{j=1}^{\infty}\zeta_j\ep_{t-j}$. If we assume, without loss of generality, that $\zeta_1\neq 0$ and also $\zeta_1=1$ for simplicity, we have 
$$f_m(z)=\int f_{\varepsilon}(z-x)\nu(dx),$$
where $\nu$ denotes the probability distribution of the random variable $\zeta_0+\sum_{j=2}^{\infty}\zeta_j\varepsilon_{t-j}$. Hence the boundedness and Lipschitz property of $f_m$ follows from assumption ${\bf A5}$. 
\item
Finally, assumptions ${\bf A4,A8}$ hold true using for instance conditional maximum likelihood estimators. See \citet{BD}, Chapter $8$, for some asymptotic results for different inference methods of ARMA parameters. 
\end{itemize}
Let us now compare our results with that of \citet{SW2}. The results obtained by these authors are very general for applications to linear processes which contain ARMA processes as a special case. Their results are shaper than ours because they obtained uniform convergence of their convolution estimate on the real line whereas we consider uniformity only on compact intervals. However, our results applies if $\varepsilon_0$ has a moment of order $3$, whereas \citet{SW2} use a moment of order $4$ (see assumption F of the paper). Moreover, the kernels used in \citet{SW2} cannot be nonnegative (see the assumption $K$ applied with an order $m\geq 2$ for the kernel) and then the estimator of the density can take negative values. This excludes some classical kernels often used by the practitioners.

\subsection{Pure GARCH models}
In this subsection, we consider the process 
$$X_t=m_0+\ep_t\sigma_t(\theta_0),\quad \sigma_t^2(\theta_0)=\alpha_{00}+\sum_{j=1}^Q\alpha_{0j}(X_{t-j}-m_0)^2+\sum_{j=1}^P\beta_{0j}\sigma^2_{t-j}(\theta_0),$$
with $\E\ep_0=0$, $\E\ep_0^2=1$.
We set $\theta=\left(m,\alpha_0,\ldots,\alpha_Q,\beta_1,\ldots,\beta_P\right)$. Moreover let 
$$\sigma^2_t(\theta)=\alpha_0+\sum_{j=1}^Q\alpha_j(X_{t-j}-m)^2+\sum_{j=1}^Q\beta_j\sigma^2_t(\theta).$$
Then $(X_t)_{t\in\Z}$ is (up to parameter $m_0$) a GARCH$(p,q)$ process.
We set
$$\underline{Y}_t=\left((X_t-m_0)^2,\ldots,(X_{t-Q+1}-m_0)^2,\sigma^2_t(\theta_0),\ldots,\sigma^2_{t-P+1}(\theta_0)\right)'.$$
There exist a sequence of i.i.d random matrices $\left(A_t\right)_{t\in\Z}$ of size $(p+q)\times (p+q)$ and a sequence of random vectors $(B_t)_{t\in\Z}$ of dimension $p+q$ such that 
\begin{equation}\label{iterG}
\underline{Y}_t=A_t\underline{Y}_{t-1}+B_t,\quad t\in\Z.
\end{equation}
We defer the reader to \cite{FZ}, p.$29$, for a precise expression of $(A_t,B_t)$ as well as the definition of the Lyapunov exponent $\gamma(A)$ of the sequence $(A_t)_{t\in\Z}$. The following assumptions will be needed.

\begin{description}
\item[G1] $\gamma(A)<0$ and for all $\theta\in \Theta$, $\sum_{j=1}^P\beta_j<1$.
\item[G2] We have $\alpha_{0,j}>0$ and $\beta_{0,j'}>0$ for $0\leq j\leq Q$ and $1\leq j'\leq P$. 
\end{description}
In the sequel, we denote by $C$ a generic positive constant.
Under the assumption {\bf G1}, there exist $s>0$ and an integer $k\geq 1$ such that $c=\E^{1/k}\left(\Vert A_kA_{k-1}\cdots A_1\Vert^s\right)<1$ and $\E\sigma_t^{2s}(\theta_0)<\infty$.
Using the representation 
$$\underline{Y}_t=B_t+\sum_{j=1}^{\infty}A_t\cdots A_{t-j+1}B_{t-j}$$
and the fact that $A_t,B_t\in\sigma(\ep_t)$, we get
\begin{eqnarray*}
\E\Vert \underline{Y}_t-\underline{Y}_{t\ell}\Vert^s&\leq& 2\sum_{j\geq \ell}\E\Vert A_t\cdots A_{t-j+1}B_{t-j}\Vert^s\\
&\leq& C\sum_{j\geq \ell}c^j\\
&\leq& C c^{\ell}. 
\end{eqnarray*}
Then, we have also $\E\vert X^2_t-X^2_{t\ell}\vert^s\leq C c^{\ell}$.
Now we check the assumptions {\bf A2}, {\bf A3}, {\bf A4} and {\bf A8}, {\bf A5} and {\bf A6}.
\begin{enumerate}
\item
We first check {\bf A2}. Setting for $t\in\Z$,
$$\underline{\sigma}^2_t(\theta)=\left(\sigma_t^2(\theta),\sigma^2_{t-1}(\theta),\ldots,\sigma^2_{t-P+1}(\theta)\right)^T$$
we have the recursive equations $\underline{\sigma}^2_t(\theta)=A_2(\theta)\underline{\sigma}_t^2(\theta)+B_{2,t}(\theta)$ with
$$B_{2,t}(\theta)=\left(\alpha_0+\sum_{j=1}^Q\alpha_j(X_{t-j}-m)^2,0,\ldots,0\right)^T,\quad A_2(\theta)=\begin{pmatrix}\beta_1&\cdots& \beta_P\\ &&0\\&I_{P-1}&\vdots\\&&0\end{pmatrix}.$$
Then $A_2(\theta)$ is the companion matrix associated to $\mathcal{P}(z)=1-\sum_{j=1}^P\beta_jz^j$.
Since $\sum_{j=1}^P\beta_{0j}<1$, the spectral radius of $A_2(\theta)$ is less than $1$ and there exists a positive integer $k$ such that $\Vert A_2(\theta_0)^k\Vert<1$ and then $\rho^k=\nor{A^k_2}=\sup_{\theta\in \Theta_{0,\epsilon}}\Vert A_2(\theta)^k\Vert<1$ if $\epsilon>0$ is small enough, using the continuity of $\theta\mapsto A_2(\theta)$. 
In particular, we have the expansion
$$\underline{\sigma}^2_t(\theta)=\sum_{j\geq 0}A_2(\theta)^jB_{2,t-j}(\theta).$$
Then, using the fact that the sequence $\left(\nor{A_2^j}\right)_{j\geq 1}$ is bounded and $\E X_t^{2s}<\infty$, we deduce that 
\begin{eqnarray*}
\E\sup_{\theta\in\Theta_{0,\epsilon}}\Vert\underline{\sigma}^2_t(\theta)-\underline{\sigma}^2_{t\ell}(\theta)\Vert^s
&\leq& C\left[\sum_{j=1}^{\frac{\ell}{2}}c^{(\ell-j)s}+\sum_{j>\frac{\ell}{2}+1}\nor{A_2^j}^s\right]\\
&\leq&  C\left(c^{\frac{\ell s}{2}}+\rho^{\frac{s\ell}{2}}\right)\\
&\leq& C \left[c^{\frac{s}{2}}\vee \rho^{\frac{s}{2}}\right]^{\ell}.
\end{eqnarray*}
This shows that $\E\left[\max_{\theta\in\Theta_{0,\epsilon}}\left\vert \sigma_t^2(\theta)-\sigma_{t\ell}^2(\theta)\right\vert^s\right]\leq C a^{\ell}$ for some $a\in(0,1)$. The proof of 
$$\E\left[\max_{\theta\in\Theta_{0,\epsilon}}\sigma_t^{2s}(\theta)\right]<\infty,\quad \E\left[\max_{\theta\in\Theta_{0,\epsilon}}\left\vert\sigma_t^2(\theta)-\overline{\sigma}_t^2(\theta)\right\vert^s\right]<\infty$$ 
is similar, using the expansion of $\underline{\sigma}^2_t(\theta)$ and the fact that $\E X_t^{2s}<\infty$. 
\item 
The assumption {\bf A3} follows from the fact that the random variables $\nor{\dot{\sigma_t^2},\sigma^2_t}$ and $\nor{\ddot{\sigma_t^2},\sigma^2_t}$ have moments of any order if $\epsilon$ is sufficiently small (see the proof Theorem $7.2$ in \cite{FZ}, part $c$ for the main arguments used for showing these properties). 
\item
For checking {\bf A4} and {\bf A8}, one can use the Gaussian QML estimator. When all the GARCH coefficients are assumed to be positive, the 
representation {\bf A8} holds for the corresponding estimator (see \cite{FZ}, p. $159-160$). Note that the representation {\bf A8} requires the assumptions {\bf G1} and {\bf G2}.
\item
Finally we check the assumptions {\bf A5} or {\bf A6}. 
First, we note that assumption {\bf A5} does not hold when $v=m_0$. In this case, the ratio $\frac{v-m_0}{\sigma_t(\theta_0)}$
is degenerate and does not have a density. For estimating $f(m_0)$, one can use the classical kernel estimate, the approach proposed in this paper has no interest because the convergence rate will be similar. In the sequel, we assume that $m_0\notin I$. 
Using Lemma \ref{regule} and the representation of GARCH processes as 
an ARCH$(\infty)$ process, one can see that $\E\sigma_t(\theta_0)<\infty$ is sufficient for {\bf A5}. Moreover, using the additional assumptions $\E\sigma_t(\theta_0)^e<\infty$ and $u\mapsto \vert u\vert^e f_{\ep}(u)$ is bounded for $e>\frac{3}{2}$, assumption {\bf A6}
also holds if $I\subset (m_0,\infty)$ or $I\subset (-\infty,m_0)$. Note that the moment condition $\E\sigma_t(\theta_0)^e<\infty$ is satisfied under the classical condition $\sum_{j=1}^Q\alpha_{0,j}+\sum_{j=1}^P\beta_{0,j}<1$ which implies $\E\sigma_t(\theta_0)^2<\infty$.   
\end{enumerate}

\paragraph{Notes}
\begin{enumerate}
\item
If we assume that $m=0$ in the model, one can use the logarithm to get
$$\log\left(X_t^2\right)=\log\left(\sigma_t^2(\theta_0)\right)+\log\left(\ep_t^2\right).$$
One can apply our results for estimating directly the density of $\log\left(X_t^2\right)$.
Setting $Z_t=\log\left(\ep_t^2\right)$ and $m_t(\theta_0)=\log\left(\sigma_t^2(\theta_0)\right)$,
one can show that the density of $Z_t$ satisfies the assumption {\bf A5} if we additionally assume that $u\mapsto u^3f^{''}_{\ep}(u)$ is bounded.
Moreover, we have $g_v(x)=f_m(v-x)=e^{v-x}f_{\sigma^2}\left(e^{v-x}\right)$ where $f_m$ denotes the density of $m_t(\theta_0)$. Then one can show that assumption {\bf A6} is satisfied if $\E\sigma_t^2(\theta)<\infty$, which is a classical condition found in practice in using GARCH models. 
Moreover, it is also possible to show that assumption {\bf A7} is satisfied under the additional condition: $u\mapsto u^{\frac{3}{2}+\delta}f_{\ep}(u)$ is bounded. The proof is omitted since one can use Lemma \ref{regule} (3.) as well as some arguments  
used in the proof of Lemma \ref{regule}. All the other assumptions are automatically satisfied it {\bf G1} and {\bf G2} hold true.
Note that the root $n$ consistent estimation of the density of $\log\sigma^2_t(\theta_0)$ is studied in \citet{Del}, for a GARCH$(1,1)$ process. Here we consider the estimation of the density of $\log X_t^2$ which is a different problem but this convergence rate also holds. Note also that we consider a more general GARCH$(P,Q)$ model in this work.

\item
One can prove similar results for pure ARCH processes (i.e $\beta_1=\cdots=\beta_p=0$ in the GARCH model), 
assuming $\alpha_1,\ldots,\alpha_q>0$. However assumption {\bf A6} (resp. assumption {\bf A7}) requires $q\geq 2$ (resp. $q\geq 3$). 
See Lemma \ref{regule} for details. Let us show that assumption {\bf A6} is not satisfied in the case $q=1$.
In the case $q=1$, we have $\sigma^2_t(\theta_0)=\alpha_{00}+\alpha_{01}X_{t-1}^2$. 
Then if $f_X$ is the marginal density of the process, we have 
$$f_{\sigma^2}(y)=\mathds{1}_{(\alpha_{00},\infty)}(y)\frac{\sqrt{\alpha_{01}}}{2\sqrt{y-\alpha_{00}}}f_X\left(\sqrt{\frac{y-\alpha_{00}}{\alpha{01}}}\right).$$
Let us assume that $0<\alpha_{01}<1$, $v$ is positive and $f(0)>0$ (the last condition holds true when $f_{\ep}(0)>0$). 
Then we have $\int g_v(x)^2f_{\ep}(x)dx=\infty$ when $f_{\ep}\left(\frac{v}{\sqrt{\alpha_{00}}}\right)>0$.
However, one can check that $\int g_v(x)^{2-\delta} dx<\infty$ for any $\delta\in (0,1)$. Indeed, we have
$$\int g_v^{2-\delta}(x)dx\leq C\int_{\alpha_{00}}^{\infty}\sqrt{y}f_{\sigma^2}(y)^{2-\delta}dy$$
which is finite (the integrability holds around the singularity $y=\alpha_{00}$, $f_{\sigma^2}$ is bounded outside a neighborhood of $\alpha_{0,0}$ and $\int \sqrt{y}f_{\sigma^2}(y)dy=\E\left[\sigma_t(\theta_0)\right]<\infty$).
We recover a phenomenon described by \citet{SWMetrika} when Frees estimator is applied for estimating the density of a sum of powers of two independent random variables. In general, a slower convergence rate is obtained when square integrability of the density fails.
See \citet{SWMetrika}, Theorem $2$, where the rate $\frac{n}{\log(n)}$ is obtained in the estimation of the density of a sum of squares $X_1^2+X_2^2$.   
   
\item    
When some parameters of the GARCH process are equal to zero, assumption {\bf A3} is not always guaranteed unless assuming 
$\E\left(X_t^6\right)<\infty$. In this case assumption {\bf A3} is automatic.    
\end{enumerate}

\subsection{ARMA processes with GARCH noises}
In this subsection, we consider the model
$$X_t-\eta_0=\sum_{j=1}^p a_{0j}\left(X_{t-j}-\eta_0\right)+Z_t-\sum_{j=1}^q b_{0j}Z_{t-j},\quad Z_t=\ep_t\sigma_t(\theta_0),$$
$$\sigma_t^2(\theta_0)=\alpha_{00}+\sum_{j=1}^Q\alpha_{0j} Z_{t-j}^2+\sum_{j=1}^P\beta_{0j}\sigma_{t-j}^2(\theta_0).$$
We define for $\theta=\left(\eta,a_1,\ldots,a_p,b_1,\ldots,b_q,\alpha_0,\ldots,\alpha_Q,\beta_1,\ldots,\beta_P\right)$, 
As for ARMA processes, we define
$$Z_t(\theta)=X_t-\eta-\sum_{j=1}^p a_j \left(X_{t-j}-\eta\right)+\sum_{j=1}^qb_j Z_{t-j}(\theta),
\quad \sigma_t^2(\theta)=\alpha_0+\sum_{j=1}^Q\alpha_j Z^2_{t-j}(\theta)+\sum_{j=1}^P\beta_j \sigma^2_{t-j}(\theta).$$
In addition to assumption {\bf G1} for the Garch parameters $(\alpha_0,\ldots,\alpha_Q,\beta_1,\ldots,\beta_P)$, we consider the following classical assumption which guarantees causality and invertibility of the ARMA part.
\begin{description}
\item[AG1] 
The roots of the two polynomials $\mathcal{P}$ and $\mathcal{Q}$ defined by 
$$\mathcal{P}(z)=1-\sum_{j=1}^pa_{0,j} z^j,\quad \mathcal{Q}(z)=1-\sum_{j=1}^qb_{0,j}z^j$$
are outside the unit disc.
\item[AG2] We have $\E Z_t^6<\infty$.
\end{description}
Note the if $\epsilon>0$ is small enough, the assumption {\bf AG1} will be also valid for all $\theta\in \Theta_{0,\epsilon}$. 
Assumption {\bf AG2} is restrictive but we do not find a way to avoid this moment condition for checking assumption {\bf A3}. This restriction is due to the technique used for the proof of Theorem $2$, with a control of the derivative of our estimator with respect to $\theta$ when $\theta$ is close to $\theta_0$. Note that, under the assumptions {\bf AG1-AG2}, we have $\E X_t^6<\infty$. 

We now check the assumptions ${\bf A2-A8}$, except {\bf A5} and {\bf A7} which we were able to check only in the pure GARCH case.

\begin{enumerate}
\item
For the assumption {\bf A2}, one can use the following expansions for $\theta\in\Theta_{0,\epsilon}$. 
If 
$$\underline{X}_t=\left(X_t,X_{t-1},\ldots,X_{t-P+1}\right)^T,\quad \underline{Z}_t(\theta)=\left(Z_t(\theta),\ldots,Z_{t-Q+1}(\theta)\right)^T,$$
$$\underline{\sigma}_t^2(\theta)=\left(\sigma^2_t(\theta),\ldots,\sigma^2_{t-q+1}(\theta)\right)^T,$$
we have 
$$\underline{X}_t=A_3(\theta)\underline{X}_{t-1}-\eta A_3(\theta)\mathds{1}+\eta\mathds{1}+B_{3,t}(\theta),\quad \underline{Z}_t(\theta)=A_1(\theta)\underline{Z}_{t-1}(\theta)+B_{1,t}(\theta),$$
$$\underline{\sigma}_t^2(\theta)=A_2(\theta)\underline{\sigma}_{t-1}^2(\theta)+B_{2,t}(\theta),$$
where $A_1(\theta)$, $A_2(\theta)$ and $A_3(\theta)$ are the companion matrices associated to $(a_1,\ldots,a_p)$, $(b_1,\ldots,b_q)$ and $(\beta_1,\ldots,\beta_P)$ respectively and
$$B_{3,t}(\theta)=\begin{pmatrix} Z_t(\theta)-\sum_{j=1}^qb_jZ_{t-j}(\theta)\\0\\\vdots\\0\end{pmatrix},\quad 
B_{1,t}(\theta)=\begin{pmatrix} X_t-\eta-\sum_{j=1}^pa_j(X_{t-j}-\eta)\\0\\\vdots\\0\end{pmatrix},$$  
$$B_{2,t}(\theta)=\begin{pmatrix} \alpha_0+\sum_{j=1}^Q\alpha_jZ^2_{t-j}(\theta)\\0\\\vdots\\0\end{pmatrix}.$$
Then assumption {\bf A2} follows from the fact that if $\epsilon>0$ is small enough, there exists three positive integers $k_1,k_2,k_3$ such that $\sup_{\theta\in\Theta_{0,\epsilon}}\Vert A_i(\theta)^{k_i}\Vert<1$ for $i=1,2,3$. The proof uses the same arguments as in the pure GARCH case and is omitted.
\item
The assumption {\bf A3} holds true if we assume ${\bf AG2}$. In this case, using the expansions 
$$\underline{\sigma}^2_t(\theta)=\sum_{j=0}^{\infty}A_2(\theta)^jB_{2,t-j}(\theta),\quad \underline{Z}_t(\theta)=\sum_{j=^0}^{\infty}A_1(\theta)^jB_{1,t-j}$$
and the equation for $\underline{X}_t$ to show that
on can show that $\nor{\sigma_t^2}$, $\nor{\dot{\sigma_t^2}}$ and $\nor{\ddot{\sigma_t^2}}$ have a moment of order $3$.
Moreover, one can show that $\nor{m_i}$, $\nor{\dot{m_i}}$ and $\nor{\ddot{m_i}}$ have a moment of order $6$. This is sufficient for checking {\bf A3}.
\item
Assumptions {\bf A4} and {\bf A8} hold true if we consider the Gaussian quasi maximum likelihood. See the proof of Theorem $7.5$ in \citet{FZ}. Note that the expansion given in {\bf A8} requires that $\theta_0$ lies in the interior of $\Theta$ and that $\E Z_t^4<\infty$.  
\item
Assumption {\bf A6} is a consequence of Lemma \ref{ARMAGARCH}. Indeed, if $f_{\varepsilon}$ is bounded, the conditional density of 
$Z_t\vert Z_{t-1},Z_{t-2}\ldots$ is bounded. Since the couple $\left(m_t(\theta_0),\sigma_t(\theta_0)\right)$ can be expressed as 
\\$\left(\sum_{j=1}^{\infty}\psi_j Z_{t-j},\sqrt{\alpha_0+\sum_{j=1}^{\infty}\alpha_j Z_{t-j}^2}\right)$ for some summable sequences of coefficients $(\psi_j)_{j\geq 1}$ and $(\alpha_j)_{j\geq 1}$, Lemma \ref{ARMAGARCH} guarantees that the density $f_{m,\sigma}$ of this couple can be bounded as follows: $f_{m,\sigma}(x,y)\leq C y$ for a positive constant $C$. Then one can conclude using inequality (\ref{maborne}) and assumption {\bf AG2} which entails the condition $\E\sigma_t(\theta_0)^2<\infty$.  

\end{enumerate}

\section{Simulation study}
In this section, we compare by simulation the mean square error of our estimator with that of the classical kernel density estimate. 
Our estimator is implemented using the quadratic kernel. The standard kernel density estimate is computed using the function {\it Density} of the software R. Bandwidth selection for our estimator is beyond the scope of this paper. However, we use the simple approach proposed in \citet{Wu} which consists in multiplying the bandwidth selected for the kernel density estimate by a factor $n^{\frac{1}{5}-\kappa}$ where $\kappa$ is compatible with our theoretical results. For a bandwidth $\hat{b}=\hat{C}n^{-\frac{1}{5}}$ with the optimal rate, we then keep the constant $\hat{C}$ and simply modify the rate. In our simulations, we found that the exponent $\kappa=2/7$ provides good results.

In the sequel, we consider three simulation setups.
\begin{enumerate}
\item
In the first setup, we consider the conditionally homoscedastic case, with the AR process $X_t=0.5X_{t-1}+\varepsilon_t$ such that $(\varepsilon_t)_{t\in\Z}$ is i.i.d with $\varepsilon_1\sim t(5)$ (the Student distribution with $5$ degrees of freedom).
\item
In the second setup, we consider the pure ARCH case with a GARCH$(1,1)$ process $X_t=\varepsilon_t\sigma_t$ such that $\sigma_t^2=0.1+0.1X_{t-1}^2+0.8\sigma_{t-1}^2$.
The noise component $\varepsilon$ still follows a $t(5)$ distribution.
\item
In the last setup, we consider an AR process with a GARCH$(1,1)$ noise, 
$$X_t=0.5 X_{t-1}+Z_t,\quad Z_t=\varepsilon_t\sigma_t,\quad \sigma^2_t=0.1+0.1X_{t-1}^2+0.8\sigma_t^2.$$ 
We assume that $\varepsilon_t$ follows a standard Gaussian. One can show that the moment condition $\E Z_t^6<\infty$ required for applying our results is not satisfied in this example.
\end{enumerate}
Note that GARCH parameters are chosen so that the expectation of the square equals to $1$ and lag coefficients have typical values encountered in practice.

The marginal density is evaluated at $10$ points equally spaced, starting from $v=0$ to $v=5$.  
The true density is approximated by
$$f_X(v)\approx \frac{1}{N}\sum_{i=1}^N\frac{1}{\sigma_i}f_{\varepsilon}\left(\frac{v-m_i}{\sigma_i}\right),\quad N=500000.$$
The RMSE of the estimator is normalized by the value of the density: $\sqrt{\E\left[\hat{f}_X(v)-f_X(v)\right]^2}/f_X(v)$. This RMSE is approximated by its empirical counterpart using $10^3$ samples. GARCH parameters are estimated using the function {\it garchFit} of the package {\it fGarch}. The R code is available on request from the author.
 
\begin{figure}[H]
\includegraphics[width=5cm,height=7cm]{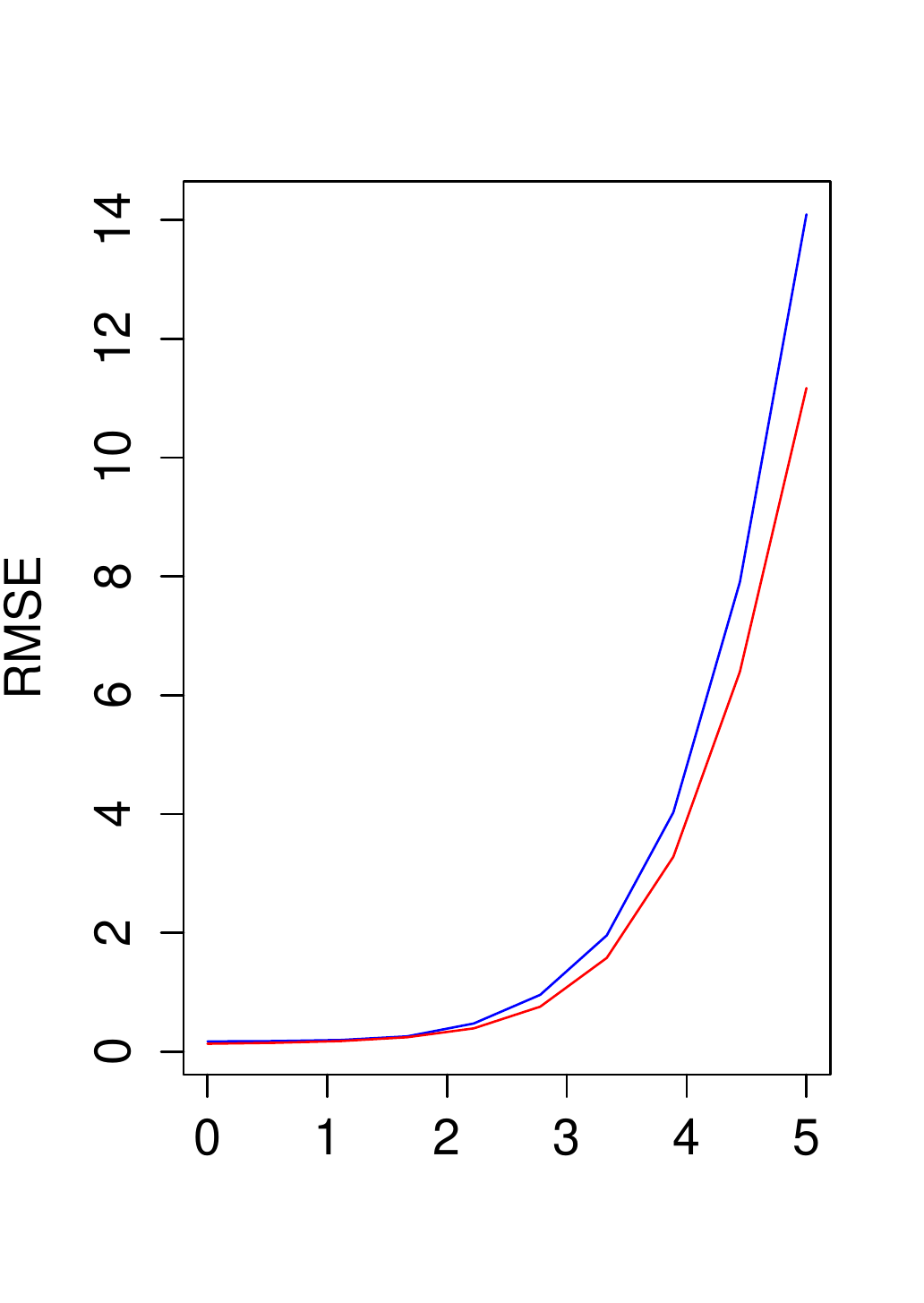}
\includegraphics[width=5cm,height=7cm]{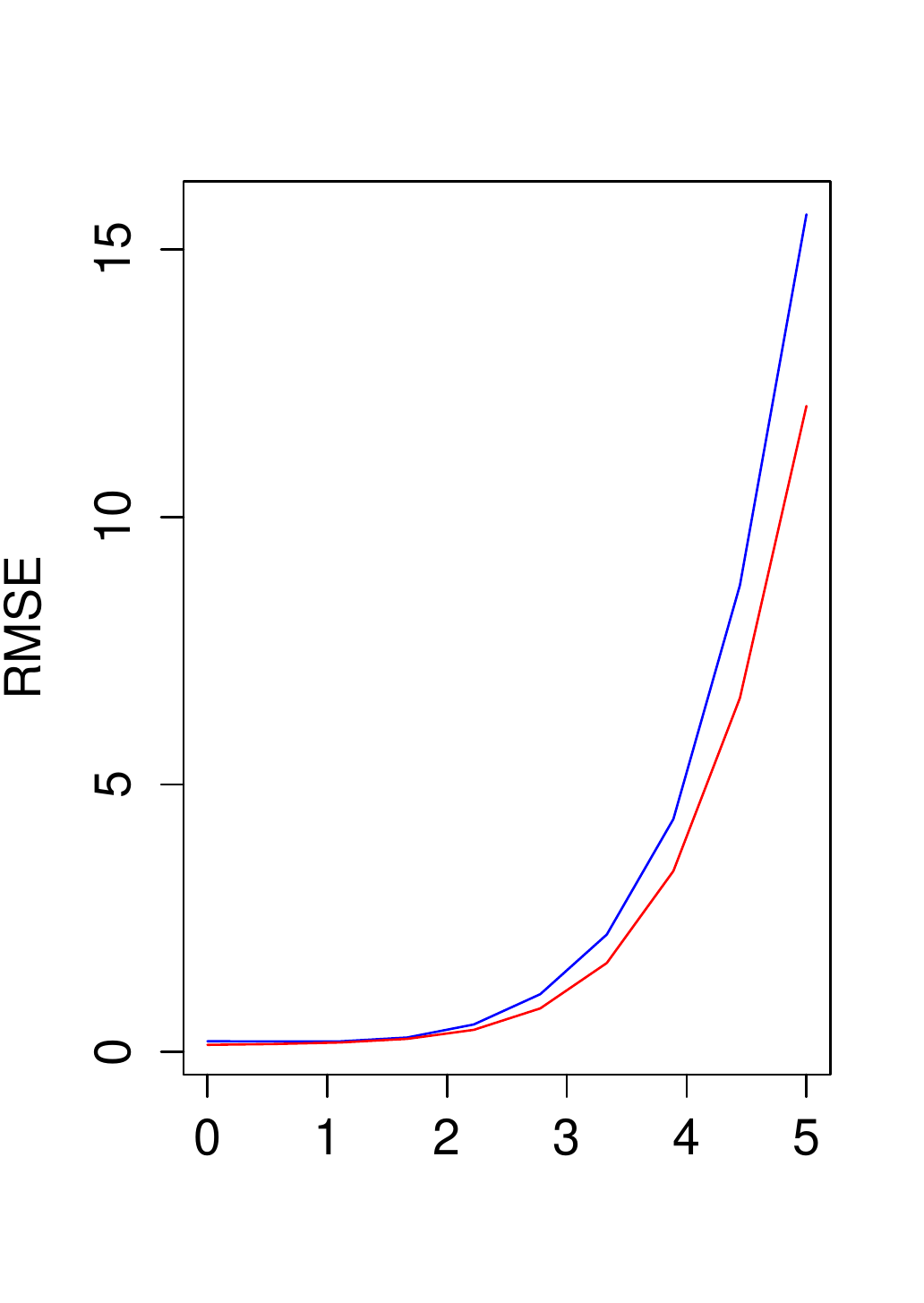}
\includegraphics[width=5cm,height=7cm]{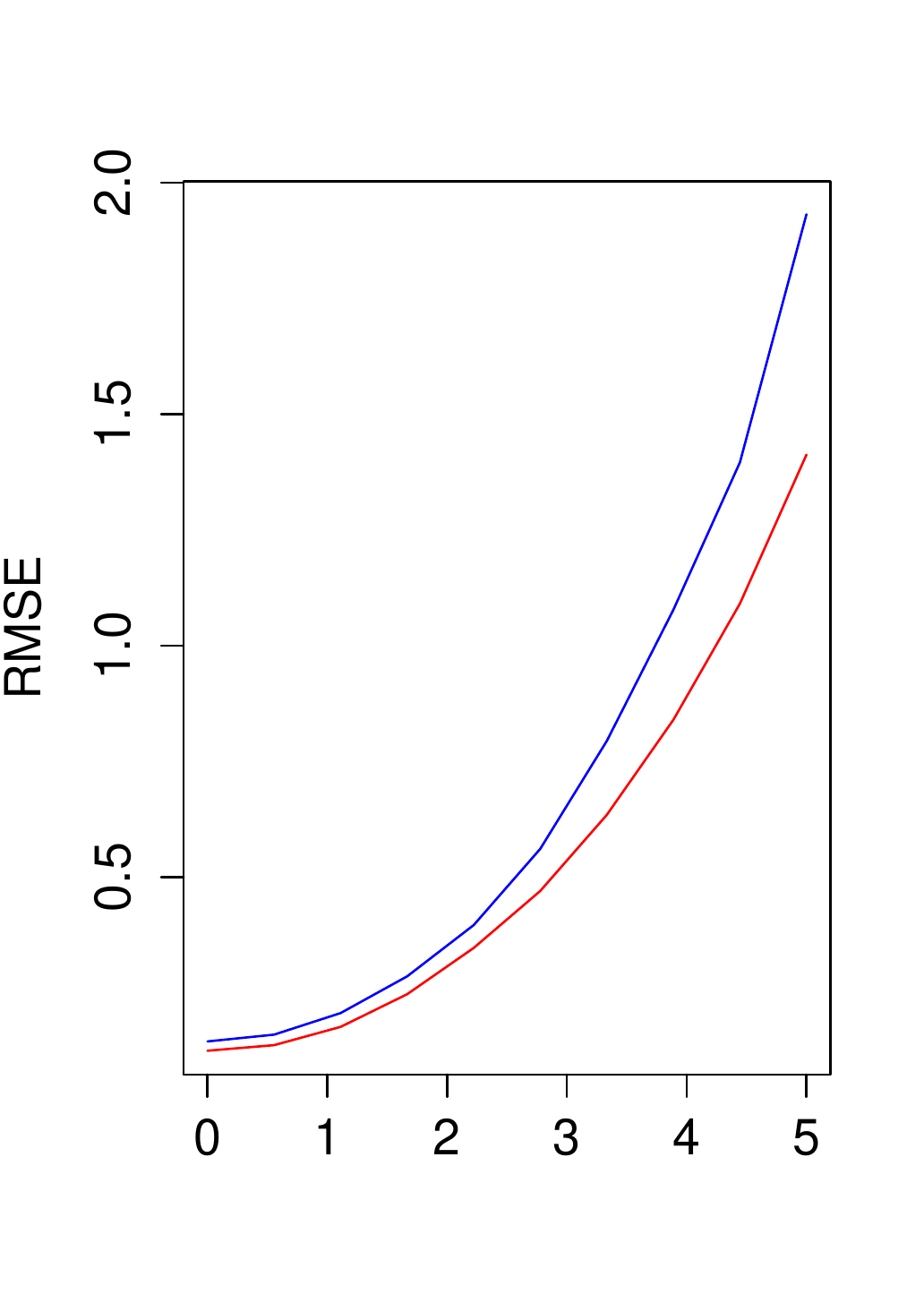}
\caption{RMSE for the AR process. Bandwidth for kernel estimate is obtained by Silverman's rule of thumb (left, $n=100$), cross-validation (middle, $n=100$) and optimal rate with constant $1$ (right, $n=100$)\label{AR}}
\end{figure}

\begin{figure}[H]
\includegraphics[width=5cm,height=7cm]{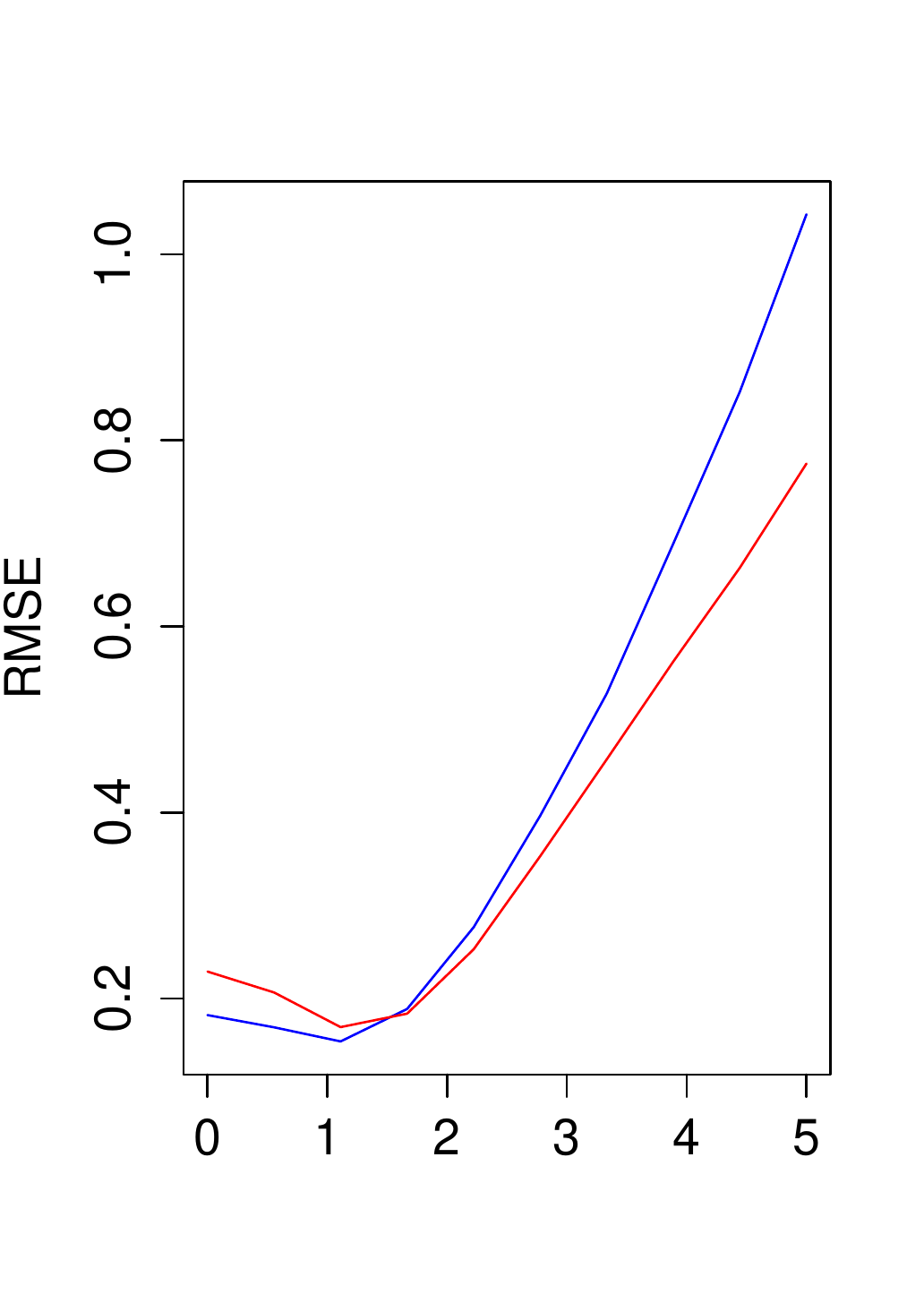}
\includegraphics[width=5cm,height=7cm]{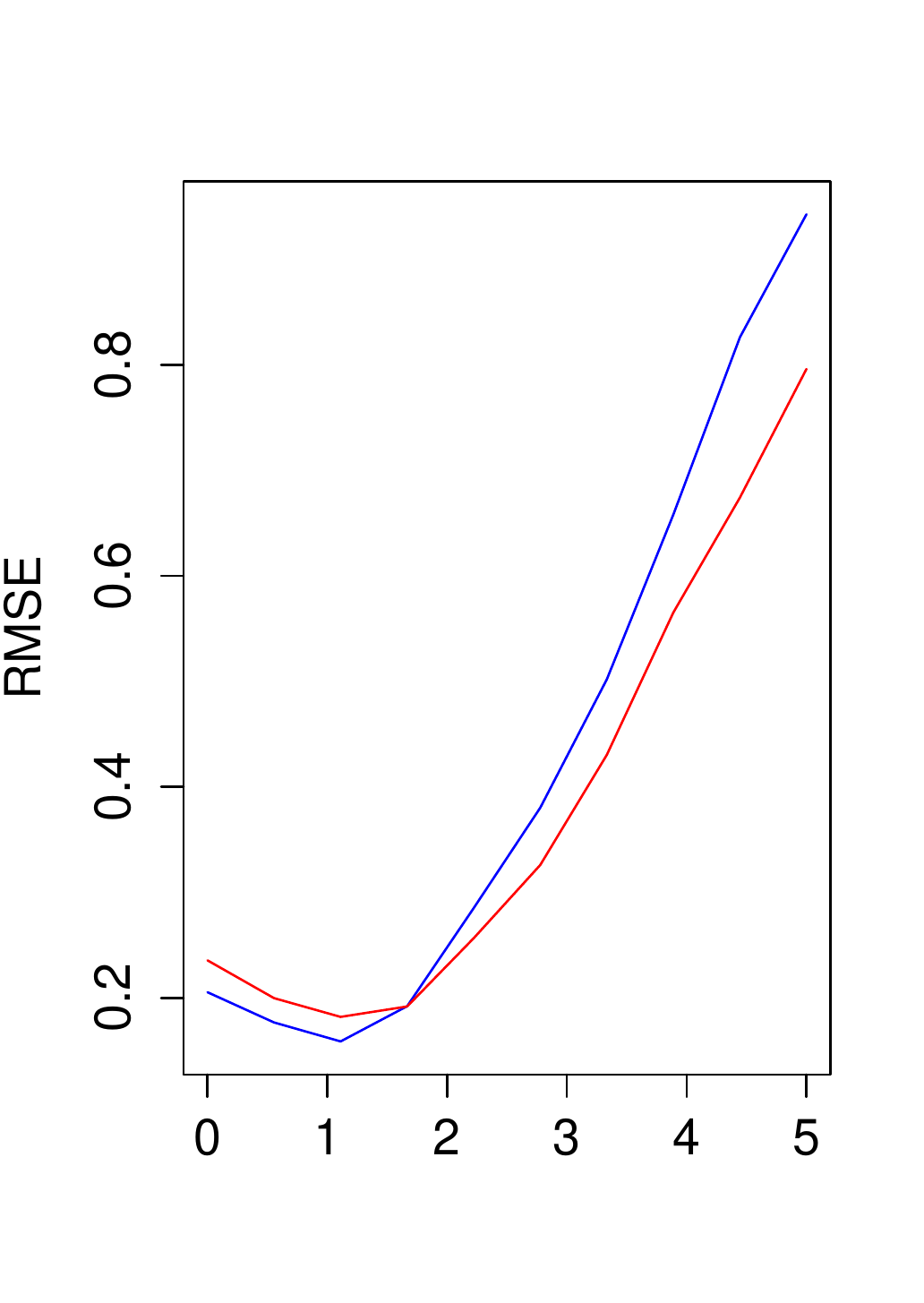}
\includegraphics[width=5cm,height=7cm]{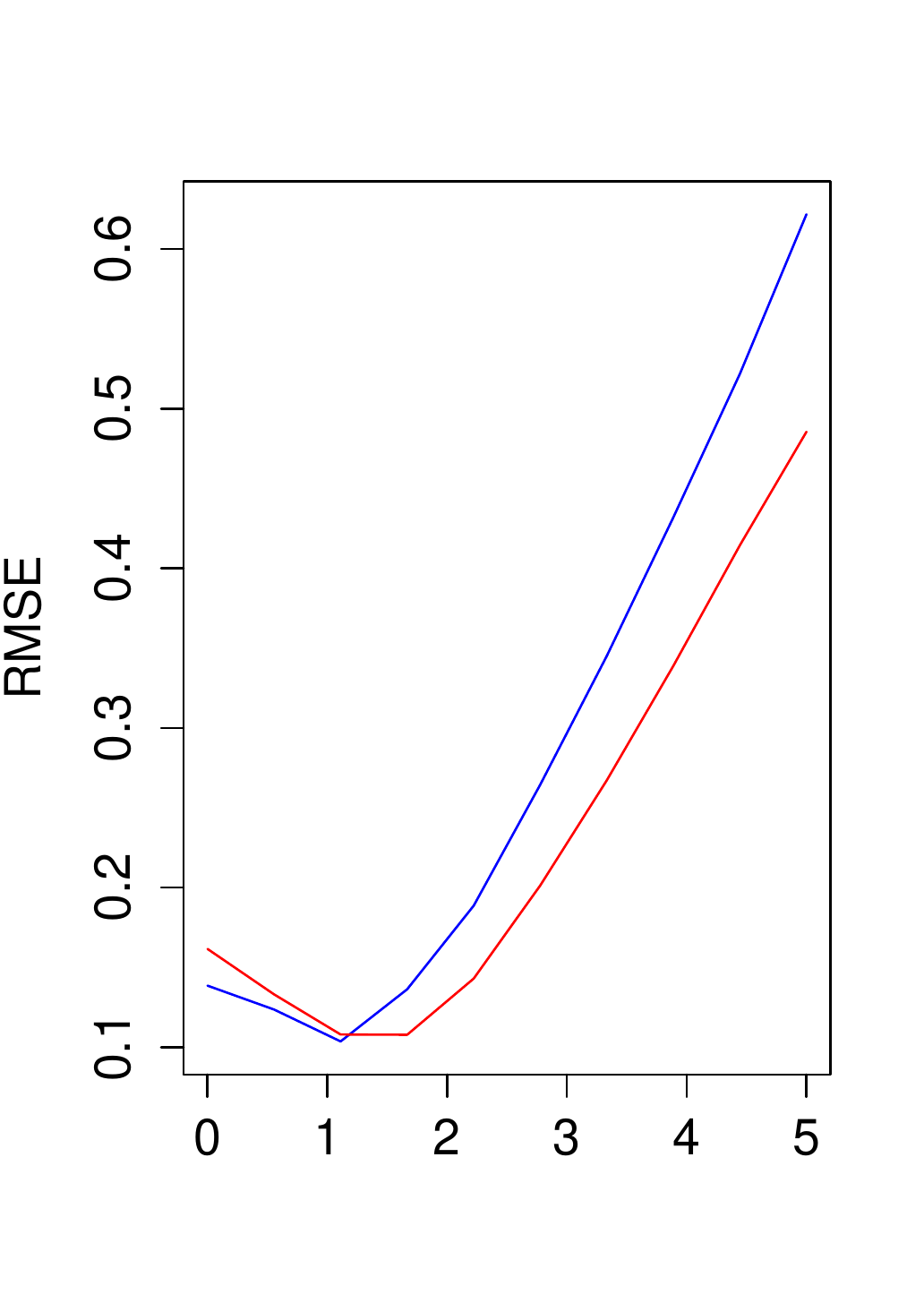}
\caption{RMSE for the GARCH process. Bandwidth is obtained using optimal rate with constant $1$ (left, $n=200$), cross-validation (middle, $n=200$) and cross-validation (right, $n=500$)\label{GARCH}}
\end{figure}

\begin{figure}[H]
\begin{center}
\includegraphics[width=5cm,height=7cm]{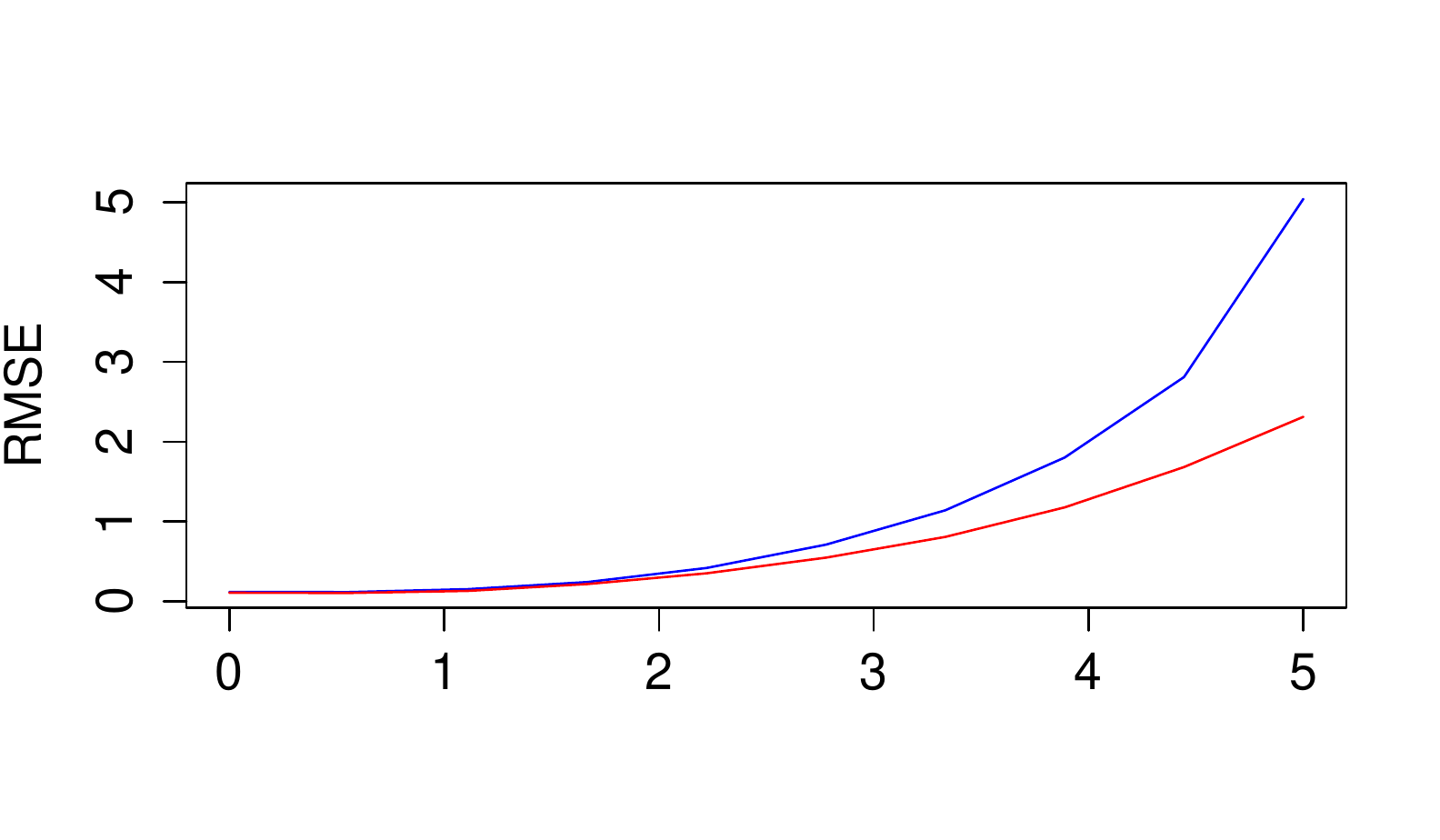}
\includegraphics[width=5cm,height=7cm]{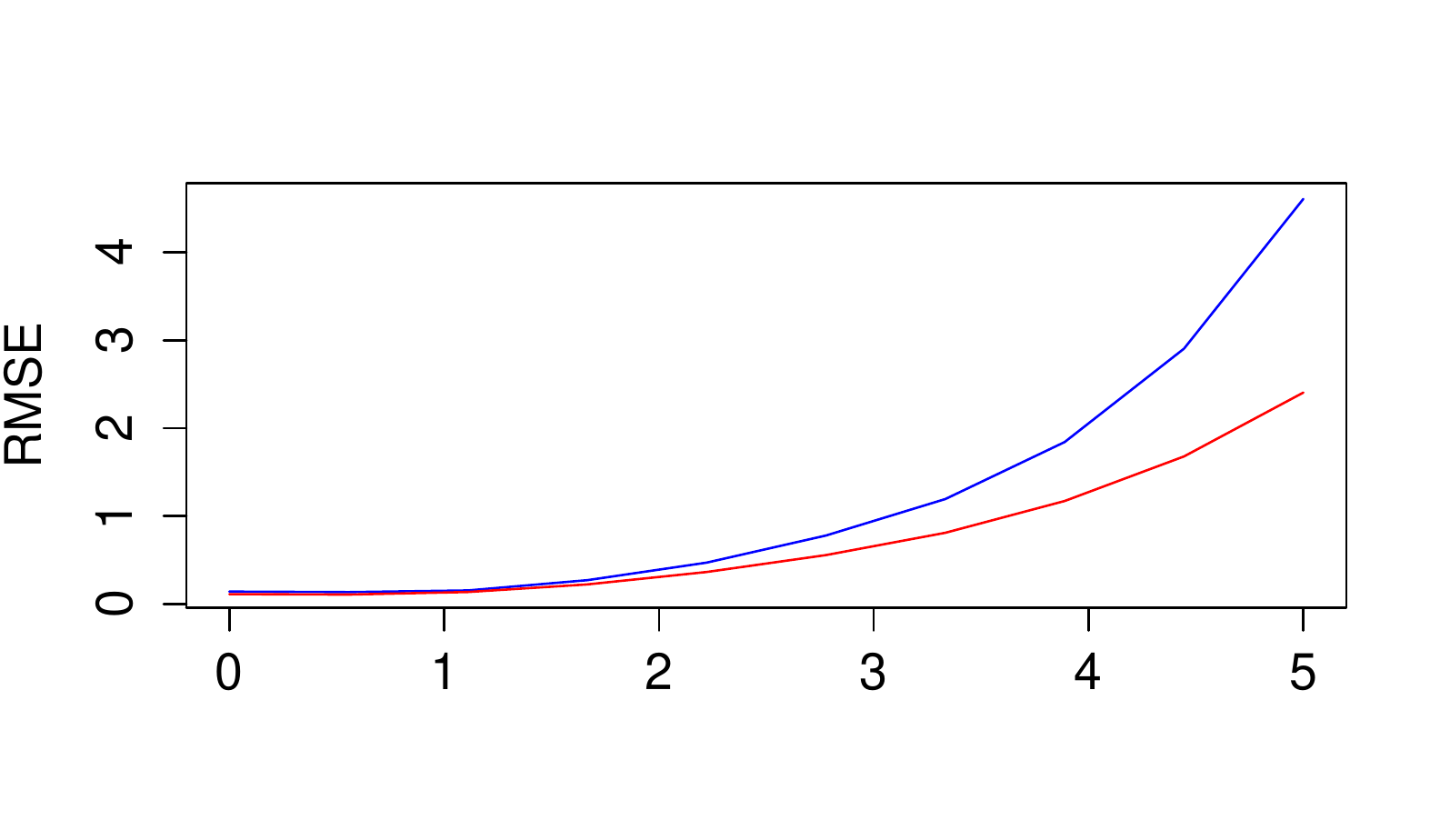}
\end{center}
\caption{RMSE for the AR-GARCH process. Bandwidth is obtained using optimal rate with constant $1$ (left, $n=200$) or cross-validation (right, $n=200$)\label{ARGARCH}}
\end{figure}

In Figure \ref{AR}, \ref{GARCH} and \ref{ARGARCH}, the blue curve represents the normalized RMSE for kernel density estimation and 
the red curve that for our method. Whatever the original bandwidth selection, our estimator performs better even if the sample size $n$ is small for estimating accurately GARCH processes. A notable exception is the second setup when $v$ is in a neighborhood of $0$. In this case, the standard method performs better. This is not surprising because of the singularity at point $v=0$, a point for which our method is almost equivalent to the standard one and our bandwidth parameter does not have the optimal convergence rate. This problem is less perceptible for the larger sample size $n=500$. A general finding is the notable superiority of our method for estimating the tails, which have an important rule in financial time series. For practical applications, adequation tests based on the marginal density as proposed in \citet{Wu} can also be adapted to GARCH processes, using our results. This problem will not be studied is the present paper.

\section{Proofs of the results}
In the subsequent proofs, $C>0$ will denote a generic constant that can change from line to line. Moreover, if $X$ is a random variable, we set $\overline{X}=X-\E(X)$. 
 
For deriving our results, we first introduce a coupling method which will be very useful. The goal of this coupling method is to construct $\ell-$dependent random sequences which approximate some weakly dependent random sequences. Then, we show that our initial estimator is asymptotically equivalent to an estimator involving an $\ell-$dependent random sequence, provided that $\ell=\ell_n$ grows with a polynomial rate.

\subsection{From dependence to $\ell-$dependence via coupling}
Let $F:\R^{\N}\rightarrow G$ be a measurable application taking values to an arbitrary measurable space $(G,\mathcal{G})$. If $Z_i=F\left(\ep_i,\ep_{i-1},\ldots\right)$,
we set 
$$Z_{i\ell}=F\left(\ep_i,\ldots,\ep_{i-\ell+1},\ep^{(i)}_{i-\ell},\ep^{(i)}_{i-\ell-1},\ldots\right),$$
where $\left\{\ep^{(i)}_t: (i,t)\in\Z^2\right\}$ is a family of i.i.d random variables independent from $(\ep_t)_{t\in\Z}$ and such that 
for all $(i,t)\in \Z^2$, $\ep^{(i)}_t$ has the same distribution than $\ep_0$.
In this case, the sequence $\left(Z_{i\ell}\right)_{i\in \Z}$ is $\ell-$dependent. This means that for all $i\in\Z$, the two $\sigma-$algebra 
$\sigma\left(Z_{j\ell}: j\geq i\right)$ and $\sigma\left(Z_{j\ell}: j\leq i-m\right)$ are independent. We call this new sequence a $\ell-$dependent approximation of $(Z_i)_{i\in\Z}$. Note that $Z_{i\ell}$ has the same distribution than $Z_i$. 
Note that the two processes $\left(m_i\right)_{i\in\Z}$ and $\left(\sigma_i\right)_{i\in\Z}$ are of this form, if $G$ denotes the set of real-valued functions defined on the set $\Theta$. We will denote by $m_{i\ell}$ and $\sigma_{i\ell}$ their corresponding $\ell-$dependent approximations. These coupling versions of the conditional mean/variance of the process will be central in our proofs.
One can note that $m_{i\ell},\sigma_{i\ell}$ and their derivatives with respect to $\theta$ have the same distribution than the original quantities. 

\subsection{A martingale decomposition}\label{martii}
The control of the derivative of our estimator will be done using appropriated martingale differences. In this subsection, 
we set for $i\in\Z$, 
$$Y_{i\ell}=\left(\ep_i,\ep_{i-1}\ldots,\ep_{i-\ell+1},\ep^{(i)}_{i-\ell},\ep^{(i)}_{i-\ell-1},\ldots\right).$$ 

Let $n'=k\ell$ and $\mathcal{I}_n=\left\{1\leq i,j\leq n': i\leq j-\ell \mbox{ or } i\geq j+\ell \right\}$.
Here $k$ denotes the integer part of the ratio $n/\ell$ and $\ell\in (0,n)$ is an integer .
for $s=1,2,\ldots,\ell$, we set 
$[s]=\left\{s+g\ell: 0\leq g \leq k-1\right\}$. For each $s$, we define two filtrations. 
We set for $g=0,1,\ldots,k-1$, 
$$\mathcal{G}_{s,g}=\sigma\left(\ep_i,\ep^{(i)}: i\leq s+g\ell\right),\quad 
\overline{\mathcal{G}}_{s,g}=\sigma\left(\ep_i,\ep^{(i)}: i> n'-s-g\ell\right).$$
Now if $T_{ij}(v,\theta)$ is a random variable measurable with respect to $\sigma\left(Y_{i\ell},Y_{j\ell}\right)$,
we set 
$$M_s(T)_g(v,\theta)=\sum_{i=1}^{s+(g-1)\ell}\left[T_{i,s+g\ell}(v,\theta)-\E\left(T_{i,s+g\ell}(v,\theta)\vert \mathcal{G}_{s,g-1}\right)\right],$$
and
$$\overline{M}_s(T)_g(v,\theta)=\sum_{i=n'-s-(g-1)\ell}^{n'}\left[T_{i,n'-s-g\ell}(v,\theta)-\E\left(T_{i,n'-s-g\ell}(v,\theta)\vert \overline{\mathcal{G}}_{s,g}\right)\right].$$
Then for each $s=1,\ldots,\ell$, $\left\{\left(M_s(T)_g,\mathcal{G}_{s,g}\right): 0\leq g\leq k-1\right\}$ and 
$\left\{\left(\overline{M}_s(T)_g,\overline{\mathcal{G}}_{s,g}\right): 0\leq g\leq k-1\right\}$ are two martingales differences.
Moreover, if $\E_{Y_{i\ell}}$ denotes integration with respect to the distribution of $Y_{i\ell}$, we have 
\begin{equation}\label{decompmart}
\sum_{(i,j)\in\mathcal{I}_n}\left[T_{ij}(v,\theta)-\E_{Y_{j\ell}}\left(T_{ij}(v,\theta)\right)\right]=\sum_{s=1}^{\ell}\left[\sum_{g=0}^{k-1}M_s(T)_g(v,\theta)+\sum_{g=0}^{k-1}\overline{M}_s(T)_g(v,\theta)\right].
\end{equation}

The following lemma will be needed in the sequel for controlling triangular arrays of martingale differences.

\begin{lem}\label{martin}
Assume that for a $\delta\in (0,1)$, $nb^{3+\delta}\rightarrow \infty$.
Let $\left(\mathcal{H}_n\right)_n$ be a sequence of finite sets such that 
$\left\vert\mathcal{H}_n\right\vert=O\left(n^{\gamma}\right)$ for some $\gamma >0$.
For each $h\in\mathcal{H}_n$,  let $\left(\xi^{(h,s)}_{n,g},\mathcal{T}^{(s)}_{n,g}\right)_{0\leq g\leq k_n}$ be a martingale difference. We assume that
$$\max_{h\in\mathcal{H}_n}\sum_{s=1}^{\ell_n}\sum_{g=1}^{k_n}\E\left(\left\vert \xi^{(h,s)}_{n,g}\right\vert^2\vert \mathcal{T}^{(s)}_{n,g}\right)=O_{\P}\left(\frac{1}{nb^3}\right)$$
and that 
$$\max_{\substack{h\in\mathcal{H}_n\\ 1\leq s\leq \ell_n\\ 1\leq g\leq k_n}}\left|\xi_{n,g}^{(h,s)}\right|=O_{\P}\left(\frac{1}{n^{2/3}b^2}\right).$$
with $\ell_n=O\left(n^t\right)$, $0<t<\frac{\delta}{2(3+\delta)}$. Then we get 
the conclusion:
$$\max_{h\in\mathcal{H}_n}\left\vert \sum_{s=1}^{\ell_n}\sum_{g=1}^{k_n}\xi_{n,g}^{(h,s)}\right\vert=o_{\P}(1).$$
\end{lem}

\paragraph{Proof of Lemma \ref{martin}.}
Let $\delta_1\in (0,2\delta/3)$ such that $t<\frac{2\delta-3\delta_1}{3(3+\delta)}$. 
We set $N=N_n=n^{-2/3} b^{-2-\delta_1}$. 
For simplicity of notations, we suppress the dependence in $n$ of our quantities.
Then we set
$$\xi^{(h,s,N)}_g=(-N)\vee \xi_{g}^{(h,s)}\wedge N-\E\left((-N)\vee \xi_{g}^{(h,s)}\wedge N\vert\mathcal{T}^{(s)}_{g-1}\right).$$
We first show that 
\begin{equation}\label{tronc+}
\max_{h\in\mathcal{H}_n}\sum_{s=1}^{\ell}\sum_{g=1}^k\left|\xi_g^{(h,s)}-\xi_g^{(h,s,N)}\right|=o_{\P}(1).
\end{equation}
To this end, we use the bound
\begin{eqnarray*}
\left|\xi_g^{(h,s)}-\xi_g^{(h,s,N)}\right|&\leq& \left|\xi_g^{(h,s)}\right|\mathds{1}_{\left|\xi_g^{(h,s)}\right|>N}+
\E\left[\left|\xi_g^{(h,s)}\right|\mathds{1}_{\left|\xi_g^{(h,s)}\right|>N}\vert \mathcal{T}_{g-1}^{(s)}\right]\\
&\leq& \left|\xi_g^{(h,s)}\right|\mathds{1}_{\left|\xi_g^{(h,s)}\right|>N}+\frac{1}{N}\E\left[\left|\xi_g^{(h,s)}\right|^2\vert \mathcal{T}_{g-1}^{(s)}\right].
\end{eqnarray*} 
If $\epsilon>0$, we have 
\begin{eqnarray*}
&&\P\left(\max_{h\in \mathcal{H}}\sum_{\ell=1}^{\ell}\sum_{g=1}^k\left|\xi_g^{(h,s)}\right|\mathds{1}_{\left|\xi_g^{(h,s)}\right|>N}>\epsilon\right)\\
&\leq& \P\left(\max_{\substack{h\in\mathcal{H}\\ 1\leq s\leq \ell\\ 1\leq g\leq k}}\left|\xi^{(h,s)}_g\right|>N\right).
\end{eqnarray*}
Using the expression of $N$ and the assumptions of the lemma, the latest probability is $o_{\P}(1)$. Moreover 
$$\frac{1}{N}\max_{h\in\mathcal{H}_n}\sum_{s=1}^{\ell_n}\sum_{g=1}^{k_n}\E\left(\left\vert \xi^{(h,s)}_{n,g}\right\vert^2\vert \mathcal{T}^{(s)}_{n,g-1}\right)=O_{\P}\left(\frac{1}{Nnb^3}\right)=o_{\P}(1).$$
This show (\ref{tronc+}). To end the proof we have to show that 

\begin{equation}\label{troncexp}
\max_{h\in\mathcal{H}}\sum_{s=1}^{\ell}\left|\sum_{g=1}^k\xi_g^{(h,s,N)}\right|=o_{\P}(1).
\end{equation}
We will prove (\ref{troncexp}) using the exponential inequality of Freedman for martingales (see \citet{Fr}). Let $\epsilon,\epsilon'>0$. First
we choose $M>0$ such that 
$$\P\left(\max_{h\in\mathcal{H}_n}\sum_{s=1}^{\ell}\sum_{g=1}^k\E\left(\left\vert \xi^{(h,s)}_{n,g}\right\vert^2\vert \mathcal{T}^{(s)}_{n,g-1}\right)>\frac{M}{nb^3}\right)\leq \epsilon'.$$
Using the fact that 
$$\E\left(\left\vert \xi^{(h,s,N)}_{g}\right\vert^2\vert \mathcal{T}^{(s)}_{g-1}\right)\leq \E\left(\left\vert \xi^{(h,s)}_g\right\vert^2\vert \mathcal{T}^{(s)}_{g-1}\right),$$
we get 
\begin{eqnarray*}
&&\P\left(\max_{h\in\mathcal{H}}\sum_{s=1}^{\ell}\left|\sum_{g=1}^k\xi_g^{(h,s,N)}\right|>\epsilon, \max_{h\in\mathcal{H}_n}\sum_{s=1}^{\ell}\sum_{g=1}^k\E\left(\left\vert \xi^{(h,s)}_{n,g}\right\vert^2\vert \mathcal{T}^{(s)}_{g-1}\right)\leq \frac{M}{nb^3}\right)\\
&\leq& \sum_{h\in\mathcal{H}}\sum_{s=1}^{\ell}\P\left(\left|\sum_{g=1}^k\xi_g^{(h,s,N)}\right|>\frac{\epsilon}{\ell}, \sum_{g=1}^k\E\left[\left\vert\xi_g^{(h,s,N)}\right\vert^2\vert \mathcal{T}^{(s)}_{g-1}\right]\leq \frac{M}{nb^3}\right)\\
&\leq& \vert 2\ell\mathcal{H}\vert\exp\left(\frac{-\epsilon^2}{\frac{2M\ell^2}{nb^3}+\frac{4\epsilon\ell N}{3}}\right). 
\end{eqnarray*}
Then (\ref{troncexp}) follows from our conditions on $b,\ell,N,\mathcal{H}$. The proof of Lemma \ref{martin} is now complete.$\square$

\subsection{Proof of Theorem $1$}
For simplicity of notations, we drop the parameter $\theta_0$ and simply write for instance $\sigma_i$ instead of $\sigma_i(\theta_0)$.
Let $0<\delta_1<\delta$ and $t\in (0,1/2)$ sufficiently small such that $\frac{n^t}{nb^{2+\delta_1}}\rightarrow 0$ and $\frac{n^t}{\sqrt{n}b}\rightarrow 0$. We denote by $\ell$ the integer part of $n^t$ and by $k$ the integer part of the ratio $n/\ell$. Then we set $n'=k\ell$ and 
$$\mathcal{I}_n=\left\{1\leq i,j\leq n': i\leq j-\ell\mbox{ or } i\geq i+\ell\right\}.$$
Note that the cardinal $\left\vert \mathcal{I}_n\right\vert$ of the set $\mathcal{I}_n$ satisfies $(n-\ell)(n-3\ell-1)\leq \left\vert\mathcal{I}_n\right\vert\leq n^2$.
We set 
$$\check{f}_{\ell}(v)=\frac{1}{n^2}\sum_{(i,j)\in\mathcal{I}_n}A_{v,ij},\quad A_{v,ij}=\frac{1}{\sigma_{i\ell}}K_b\left(\frac{v-m_{i\ell}}{\sigma_{i\ell}}-\varepsilon_j\right).$$
Note first that 
\begin{equation}\label{firstapprox}
\sqrt{n}\sup_{v\in I}\left\vert\check{f}_X(v)-\check{f}_{\ell}(v)\right\vert=o_{\P}(1).
\end{equation}
Indeed, using Lemma \ref{approx2}, we have 
$$\sqrt{n}\sup_{v\in I}\left\vert \check{f}_X(v)-\frac{1}{n^2}\sum_{1\leq i,j\leq n}A_{v,ij}\right\vert=o_{\P}(1).$$
Moreover, there exists a constant $C>0$ such that
$$\sqrt{n}\sup_{v\in I}\left\vert\check{f}_{\ell}(v)-\frac{1}{n^2}\sum_{1\leq i,j\leq n}A_{v,ij}\right\vert\leq \frac{C\ell}{\sqrt{n}b}\rightarrow 0.$$
Hence (\ref{firstapprox}) follows.
In the sequel, we study the behavior of the estimator $\check{f}_{\ell}(v)$.
\subsubsection{Bias part}
We remind that $f_X(v)=\E\left[\frac{1}{\sigma_i}f_{\ep}\left(\frac{v-m_i}{\sigma_i}\right)\right]$. 
Since $f^{''}_{\ep}$ is bounded, there exists $C>0$ such that for all $x,h\in \R$, 

\begin{equation}\label{interes}
\left\vert f_{\ep}(x+h)-f_{\ep}(x)-hf'_{\ep}(x)\right\vert \leq C h^2.
\end{equation}
From (\ref{interes}), we deduce that
\begin{eqnarray*}
&&\E\left[A_{v,ij}\right]-f(v)\\&=& \E\int \frac{1}{\sigma_i}K(w)\left[f_{\ep}\left(\frac{v-m_i}{\sigma_i}-bw\right)-f_{\ep}\left(\frac{v-m_i}{\sigma_i}\right)\right]dw\\
&=& O\left(b^2\right).
\end{eqnarray*} 
Using the condition $\sqrt{n}b^2\rightarrow 0$, we get 
$$\sqrt{n}\left(\E\left[\check{f}_{\ell}(v)\right]-f(v)\right)=o(1).$$

\subsubsection{The variance part}
We first focus on 
$$\frac{1}{n^{3/2}}\sum_{j=\ell}^{n'-1}\sum_{i=n'-j+\ell}^{n'} A_{v,i(n'-j)}.$$
We set for $j\in\N$,
$$\mathcal{H}_{nj}=\sigma\left(\ep_s,\ep^{(s)}: s\geq n'-j\right).$$
For this part we set 
$$Z_i(s,v)=\frac{1}{\sigma_{i\ell}}K_b\left[\frac{v-m_{i\ell}}{\sigma_{i\ell}}-s\right],$$
$\overline{Z}_i(s,v)=Z_i(s,v)-\E\left[Z_i(s,v)\right]$ and 
\begin{eqnarray*}
S_{v,ij}&=&A_{v,i(n'-j)}-\E\left(A_{v,i(n'-j)}\big\vert \ep_{n'-j}\right)\\&-&\E\left[A_{v,i(n'-j)}-\E\left(A_{v,i(n'-j)}\big\vert \ep_{n'-j}\right)\big\vert\mathcal{H}_{n(j-1)}\right].
\end{eqnarray*}
Using the independence properties, observe that 
\begin{eqnarray*}
S_{v,ij}&=& A_{v,i(n'-j)}-\frac{1}{\sigma_{i\ell}}\int K(w)f_{\ep}\left[\frac{v-m_{i\ell}}{\sigma_{i\ell}}-bw\right]dw\\
&-& \int K(h)g_v\left(\ep_{n'-j}+bh\right)dh+\E\left[A_{v,i(n'-j)}\right].
\end{eqnarray*}

From the first writing of $S_{v,ij}$, we have for $v,\bar{v}\in I$, 
\begin{eqnarray*}
&&\E\left\vert\sum_{j=\ell}^{n'-1}\sum_{i=n'-j-\ell}^{n'}\left[S_{v,ij}-S_{\bar{v},ij}\right]\right\vert^2\\
&\leq& 2\sum_{j=\ell}^{n'-1}\int\E\left\vert \sum_{i=n'-j-\ell}^{n'}\left[Z_i(s,v)-Z_i(s,\bar{v})\right]\right\vert^2f_{\ep}(s)ds.\\
\end{eqnarray*}
Moreover, for $s\in \R$, we have, using $\ell-$dependence,
$$\E\left\vert \sum_{i=n'-j-\ell}^{n'}\left[Z_i(s,v)-Z_i(s,\bar{v})\right]\right\vert^2\leq \ell\sum_{h=1}^{\ell}\sum_{g=0}^{k-1}\E\left\vert Z_{h+g\ell}(s,v)-Z_{h+g\ell}(s,\bar{v})\right\vert^2.$$

Moreover, we have 
\begin{eqnarray*}
&&\E\int\left\vert Z_i(s,v)-Z_i(s,\bar{v})\right\vert^2f_{\ep}(s)ds\\
&\leq& 2 \E\int \left\vert \frac{1}{\sigma_{i\ell}}K_b\left[\frac{v-m_{i\ell}}{\sigma_{i\ell}}-s\right]
-\frac{1}{\sigma_{i\ell}}K_b\left[\frac{\bar{v}-m_{i\ell}}{\sigma_{i\ell}}-s\right]\right\vert^2f_{\ep}(s)ds\\
&\leq& \frac{C\vert v-\bar{v}\vert^{\delta_1}}{b^{1+\delta_1}}\int \left\vert K(z)-K\left(\frac{\bar{v}-v}{b\sigma_{i\ell}}+z\right)\right\vert^{2-\delta_1}f_{\ep}\left(\frac{v-m_{i\ell}}{\sigma_{i\ell}}-b z\right)dz\\
&\leq& \frac{C\vert v-\bar{v}\vert^{\delta_1}}{b^{1+\delta_1}}\int\int \mathds{1}_{w\in [z,z+\frac{\bar{v}-v}{b\sigma_{i\ell}}]}\vert K'(w)\vert^{2-\delta_1} dwdz\\
&\leq& \frac{C}{b^{2+\delta_1}}\vert v-\bar{v}\vert^{1+\delta_1}.
\end{eqnarray*}
Then we get, 
\begin{equation}\label{tension}
\frac{1}{n^3}\E\left\vert\sum_{j=\ell}^{n'-1}\sum_{i=n'-j-\ell}^{n'}\left[S_{v,ij}-S_{\bar{v},ij}\right]\right\vert^2
\leq \frac{C\ell}{n b^{2+\delta_1}}\vert v-\bar{v}\vert^{1+\delta_1}.
\end{equation}
Note that $\frac{\ell}{nb^{2+\delta_1}}\rightarrow 0$. Moreover, it is easily seen that
$$\frac{1}{n^3}\E\left\vert\sum_{j=\ell}^{n-1}\sum_{i=n-j-\ell}^nS_{v,ij}\right\vert^2
\leq \frac{C\ell}{n b^2}=o(1).$$
This means that $G_n:v\mapsto\frac{1}{n^{3/2}}\sum_{j=\ell}^{n-1}\sum_{i=n-j-\ell}^nS_{v,ij}$ converges pointwise to $0$ in probability. 
Note that $G_n$ is a random function taking values in $\mathcal{C}(I)$, the space of real-valued and continuous functions defined on the compact interval $I$. From (\ref{tension}) and the Kolmogorov-Chentsov tightness criterion in $\mathcal{C}(I)$ (see for instance \citet{Kal}, Corollary $14.9$), we deduce that
$$\sup_{v\in I}\frac{1}{n^{3/2}}\left\vert\sum_{j=\ell}^{n-1}\sum_{i=n-j-\ell}^nS_{v,ij}\right\vert=o_{\P}(1).$$
Using the same type of arguments, one can also show that 
$$\sup_{v\in I}\frac{1}{n^{3/2}}\left\vert \sum_{j=\ell+1}^{n'}\sum_{i=1}^{j-\ell}\bar{S}_{v,ij}\right\vert=o_{\P}(1).$$
where 
$$\bar{S}_{v,ij}=A_{v,ij}-\int\frac{K(w)}{\sigma_{i\ell}}f_{\ep}\left[\frac{v-m_{i\ell}}{\sigma_{i\ell}}-bw\right]dw-\int K(h) g_v\left(\ep_j+bh\right)dh+\E\left[A_{v,ij}\right].$$
Moreover, since $f_{\ep}$ is $\mathcal{C}^2$, we have (see also the control of the biais) 
$$\frac{1}{n^{3/2}}\sup_{v\in I}\left\vert \sum_{(i,j)\in\mathcal{I}_n}\int\frac{K(w)}{\sigma_{i\ell}}\left[f_{\ep}\left[\frac{v-m_{i\ell}}{\sigma_{i\ell}}-bw\right]-f_{\ep}\left[\frac{v-m_{i\ell}}{\sigma_{i\ell}}\right]\right]dw\right\vert=O_{\P}\left(\sqrt{n}b^2\right).$$
Now we are going to show that 
\begin{equation}\label{puisi}
\frac{1}{n^{3/2}}\sup_{v\in I}\left\vert \sum_{(i,j)\in\mathcal{I}_n}\int K(h)\left[g_v(\ep_j+bh)-g_v(\ep_j)\right]\right\vert=o_{\P}(1).
\end{equation}
The proof is divided into two parts. In the first part, we show that 
\begin{equation}\label{puisi1}
\frac{1}{n^{3/2}}\sup_{v\in I}\left\vert \sum_{(i,j)\in\mathcal{I}_n}\int K(h)\E\left[g_v(\ep_j+bh)-g_v(\ep_j)\right]dh\right\vert=o_{\P}(1).
\end{equation}
We have 
\begin{eqnarray*}
\int K(h)\E\left[g_v(\ep_j+bh)-g_v(\ep_j)\right]dh&=&\int\int K(h)\left[g_v(u+bh)-g_v(u)\right]f_{\ep}(u)dudh\\
&=& \int \int K(h) g_v(z)\left[f_{\ep}(z-bh)-f_{\ep}(z)\right]dhdz.
\end{eqnarray*}
Using (\ref{interes}), the condition $\sqrt{n}b^2\rightarrow 0$ and the fact that $\sup_{v\in I}\int g_v(z)dz\leq C \sup_{v\in I}\int f_{\frac{v-m_i}{\sigma_i}}(z)dz\leq C$, we get (\ref{puisi1}). 
To show (\ref{puisi}), it remains to show that 
\begin{equation}
\begin{split}
&\frac{1}{n^{3/2}}\sup_{v\in I}\left\vert \sum_{(i,j)\in\mathcal{I}_n}\int K(h)\left[\overline{g_v(\ep_j+bh)}-\overline{g_v(\ep_j)}\right]dh\right\vert\\
&=\frac{1}{\sqrt{n}}\sup_{v\in I}\left\vert \sum_{j=1}^{n'}c_{n,j}\int K(h)\left[\overline{g_v(\ep_j+bh)}-\overline{g_v(\ep_j)}\right]dh\right\vert\\
&= o_{\P}(1),
\end{split}
\label{puisi2}
\end{equation}
where $c_{n,j}=\frac{1}{n}\sum_{i=1}^{n'}\mathds{1}_{\vert i-j\vert\geq \ell}$.
\begin{enumerate}
\item
To show (\ref{puisi2}) when $I$ is the singleton $\{v\}$, we use Jensen inequality and the fact that the translations are continuous in $\L^2$. More precisely,
we have
\begin{eqnarray*}
&&n^{-3}\E\left\vert \sum_{(i,j)\in \mathcal{I}_n}\int K(h)\left[\overline{g_v(\ep_j+bh)}-\overline{g_v(\ep_j)}\right]dh\right\vert^2\\ 
&\leq& n^{-1}\sum_{j=1}^{n'}c_{n,j}^2\int\left\vert\int K(h)\left[g_v(u+bh)-g_v(u)\right]dh\right\vert^2f_{\ep}(u)du\\
&\leq& n'/n\int \int K(h)\left[g_v(u+bh)-g_v(u)\right]^2f_{\ep}(u)du dh\\
&\leq& \sup_{\vert z\vert\leq b}\int \left[g_v(u+z)-g_v(u)\right]^2f_{\ep}(u)du\\
&=& o_{\P}\left(1\right).
\end{eqnarray*}
\item
To show (\ref{puisi2}) when $I$ is a compact interval not reduced to one point, we proceed in two parts.
\begin{enumerate}
\item
We first show that 
\begin{equation}\label{raj}
\frac{1}{\sqrt{n}}\sup_{v\in I}\left\vert \sum_{j=1}^{n'}\int K(h)\left[g_v\left(\ep_j+bh\right)-g_v(\ep_j)\right]dh\right\vert=o_{\P}(1).
\end{equation}
To this end, we will use Lemma $19.34$ in \cite{VW}.
We set 
$$g_{n,v}(x)=\int K(h)\left[g_v(x+bh)-g_v(x)\right]dh.$$
For the family $\mathcal{G}_{n,I}=\left\{g_{n,v}:v\in I\right\}$ of functions, we consider the envelope function $G_n$ defined by 
$$G_n(x)=\int K(h)\left[G(x+bh)+G(x)\right]dh,$$
where $G$ is defined in assumption ${\bf A7}$.
For bounding the bracketing numbers of this family, we first observe that if $[f_1,f_2]$ is an $\epsilon-$bracketing in $\mathcal{G}_I$, then 
$f_1\leq g_v\leq f_2$ (i.e $\int \left(f_2-f_1\right)^2d\mu<\epsilon^2$) entails that 
$$f_{n,1}(x)=\int K(h)\left[f_1(x+bh)-f_2(x)\right]dh\leq g_{n,v}(x)\leq \int K(h)\left[f_2(x+bh)-f_1(x)\right]dh=f_{n,2}(x).$$
Moreover, one can show that
$$\int \left\vert f_{n,2}(x)-f_{n,1}(x)\right\vert^2f_{\ep}(x)dx\leq 2 \epsilon^2.$$
Then, we have 
$$N_{[]}\left(\sqrt{2}\epsilon,\mathcal{G}_{n,I},\L^2(\ep_0)\right)\leq N_{[]}\left(\epsilon,\mathcal{G}_I,\L^2(\mu)\right)$$
and the bracketing numbers of the family $\mathcal{G}_{n,I}$ are of polynomial decay. 
Next we show that 
\begin{equation}\label{rajout2}
\sup_{v\in I}\left\vert \int g_{n,v}(x)^2f_{\ep}(x)dx\right\vert=o_{\P}(1).
\end{equation}
To show (\ref{rajout2}), we consider for a given $\epsilon>0$, some brackets $I_1,I_2,\ldots,I_T$ that cover $\mathcal{G}_I$. 
For each integer $1\leq p\leq T$, we consider an element $g_{v_p}\in I_p$. Then, for $1\leq p\leq T$, if $I_p=\left[f_1^{(p)},f_2^{(p)}\right]$, we set $I_{n,p}=\left[f^{(p)}_{n,1},f_{n,2}^{(p)}\right]$, with
$$f^{(p)}_{n,1}(x)=\int K(h)\left[f^{(p)}_1(x+bh)-f^{(p)}_2(x)\right]dh,\quad f^{(p)}_{n,2}(x)=\int K(h)\left[f^{(p)}_2(x+bh)-f^{(p)}_1(x)\right]dh.$$
Then, if $g_v\in I_p$, we have $g_{n,v}\in I_{n,p}$ and 
\begin{eqnarray*}
\int g_{n,v}^2(x)f_{\ep}(x)dx&\leq& 2\int \left[g_{n,v}(x)-g_{n,v_p}(x)\right]^2f_{\ep}(x)dx+2\int g_{n,v_p}^2(x)f_{\ep}(x)dx \\
&\leq& 4\epsilon^2+2\max_{1\leq p\leq T} \int g_{n,v_p}^2(x)f_{\ep}(x)dx.
\end{eqnarray*}
Since for each $v\in I$, we have $\int g_{n,v}(x)^2f_{\ep}(x)dx=o(1)$, we conclude that 
$$\limsup_{n\rightarrow\infty}\int g_{n,v}^2(x)f_{\ep}(x)dx\leq 4\epsilon^2.$$
Since $\epsilon$ is arbitrary, we conclude (\ref{rajout2}).

Moreover, setting, for $\kappa>0$, $a_n(\kappa)=\kappa/Log\left[N_{[]}\left(\kappa,\mathcal{G}_{n,I},\L^2(\ep_0)\right)\right]$, we have
\begin{eqnarray*}
\sqrt{n}\E G_n(\ep_0)\mathds{1}_{G_n(\ep_0)>\sqrt{n}a_n(\kappa)}
&\leq& \frac{\sqrt{n}}{\left(\sqrt{n}a_n(\kappa)\right)^{1+o}}\E G_n(\ep_0)^{2+o}\\
&\leq& \frac{2}{n^{o/2}a_n(\kappa)^{1+o}}\int G(x)^{2+o}d\mu(x).
\end{eqnarray*}
Here $Log(x)=\log(x)\wedge 1$.
Then, from (\ref{rajout2}), one can choose $\kappa^2=\kappa_n^2\geq \int g_{n,v}(x)^2f_{\ep}(x)dx$ such that $\kappa\rightarrow 0$ and $n^{o/2}a_n(\kappa)^{1+o}\rightarrow \infty$. Hence, from Lemma $19.34$ in \cite{VW}, we deduce that 
$$\frac{1}{\sqrt{n}}\E\sup_{v\in I}\left\vert \sum_{j=1}^{n'}\overline{g_{n,v}(\ep_j)}\right\vert =o_{P}(1)$$
and hence (\ref{raj}).
\item
Finally, we have 
\begin{eqnarray*}
&&\sup_{v\in I}\left\vert n^{-3/2}\sum_{(i,j)\in\mathcal{I}_n}g_{n,v}(\ep_j)-\frac{1}{\sqrt{n}}\sum_{j=1}^{n'}g_{n,v}(\ep_j)\right\vert\\
&\leq& \frac{1}{\sqrt{n}}\sum_{j=1}^{n'}\left\vert c_{n,j}-1\right\vert \int K(h)\left[G\left(\ep_j+bh\right)+G\left(\ep_j\right)\right]dh\\
&\leq& \frac{1}{n}\sum_{j=1}^{n'}\int K(h)\left\vert G\left(\ep_j+bh\right)+G\left(\ep_j\right)\right\vert dh\cdot \frac{n-n'+2\ell-1}{\sqrt{n}}. 
\end{eqnarray*}
Using the fact that $\ell\sim n^t$ with $t<1/2$, we deduce from assumption ${\bf A7}$ that 
$$\sup_{v\in I}\left\vert n^{-3/2}\sum_{(i,j)\in\mathcal{I}_n}g_{n,v}(\ep_j)-\frac{1}{\sqrt{n}}\sum_{j=1}^{n'}g_{n,v}(\ep_j)\right\vert
=o_{\P}(1).$$
From the last convergence and from (\ref{raj}), we deduce (\ref{puisi2}).
\end{enumerate}
\end{enumerate}

\subsection{End of the proof of Theorem $1$}
Collecting the results of the two previous subsections, we have shown that 
$$\sqrt{n}\left(\check{f}_X(v)-f_X(v)\right)=\frac{1}{\sqrt{n}}\sum_{j=1}^{n'}c_{n,j}\overline{g_v\left(\ep_j\right)}+\frac{1}{\sqrt{n}}\sum_{i=1}^{n'}c_{n,i}\overline{\frac{1}{\sigma_{i\ell}}f_{\ep}\left(\frac{v-m_{i\ell}}{\sigma_{i\ell}}\right)}+o_{\P}(1),$$
uniformly on $I$ and with a uniform convergence on $I$ for the partial sum involving $g_v$ if assumption ${\bf A7}$ holds true.
Then it is straightforward to show that one can replace $c_{n,i}$ and $c_{n,j}$ by $1$ in this asymptotic expansion. Indeed, we have 
$$\sqrt{n}\left\vert c_{n,i}-1\right\vert \leq \frac{n-n'+2\ell-1}{\sqrt{n}}\rightarrow 0$$
and for the uniform convergence over $I$, one can use the bounds
$$g_v(\ep_j)\leq G(\ep_j),\quad \frac{1}{\sigma_{i\ell}}f_{\ep}\left(\frac{v-m_{i\ell}}{\sigma_{i\ell}}\right)\leq \frac{1}{\gamma}\Vert f_{\ep}\Vert_{\infty}.$$   
Using the same arguments, $n'$ can be replaced with $n$. Finally, using the arguments given in the proof of Lemma \ref{approx2}, we have 
$$\frac{1}{\sqrt{n}}\sup_{v\in I}\sum_{i=1}^n\left\vert \sigma_{i\ell}^{-1}f_{\ep}\left(\frac{v-m_{i\ell}}{\sigma_{i\ell}}\right)-\sigma_i^{-1}f_{\ep}\left(\frac{v-m_i}{\sigma_i}\right)\right\vert=o_{\P}(1).$$
The proof of the tightness of $v\mapsto \frac{1}{\sqrt{n}}\sum_{i=1}^n\sigma_i^{-1}\left[f_{\ep}\left(\frac{v-m_i}{\sigma_i}\right)-f_X(v)\right]$ will be studied in detail in the proof of Corollary $1$.$\square$

\subsection{Proof of Theorem $2$}
Using Lemma \ref{approx1} and our bandwidth conditions, we have 
$$\sup_{v\in I}\left\vert \hat{f}_X(v)-\widetilde{f}_X(v)\right\vert=o_{\P}(1/\sqrt{n}),$$
where $\widetilde{f}_X$ is defined as $\hat{f}_X$ but the quantities $\overline{m}_i$ and $\overline{\sigma}_i$  being replaced with $m_i$ and $\sigma_i$ respectively. 
This shows that possible truncations of the conditional mean/variance only using $X_1,\ldots,X_n$ is asymptotically negligible for our estimator.
Now, for $\theta\in \Theta$, we recall that $\left(m_{i\ell}(\theta)\right)_{i\in\Z}$ and $\left(\sigma_{i\ell}(\theta)\right)_{i\in\Z}$
denote the $\ell-$dependent approximations of $\left(m_i(\theta)\right)_{i\in\Z}$ and $\left(\sigma_i(\theta)\right)_{i\in\Z}$ respectively.
Then setting $X_{j\ell}=m_{j\ell}(\theta_0)+\ep_j\sigma_{j\ell}(\theta_0)$, we define
$$L_{v,i\ell}(\theta)=\frac{v-m_{i\ell}(\theta)}{\sigma_{i\ell}(\theta)},\quad \ep_{j\ell}(\theta)=\frac{X_{j\ell}-m_{j\ell}(\theta)}{\sigma_{j\ell}(\theta)}$$
and  
$$\widetilde{f}_{\ell}(v)=\frac{1}{n^2}\sum_{1\leq i,j\leq n}\frac{1}{\sigma_{i\ell}(\hat{\theta})}K_b\left[L_{v,i\ell}(\hat{\theta})-\ep_{j\ell}(\hat{\theta})\right].$$
In the rest of the proof, we fix a real number $t$ such that  $0<t<\frac{\delta}{2(3+\delta)}$ and we denote by $\ell$ the integer part of $n^t$. Using Lemma \ref{approx2}, we have 
$$\sup_{v\in I}\left\vert \widetilde{f}_X(v)-\widetilde{f}_{\ell}(v)\right\vert=o_{\P}\left(1/\sqrt{n}\right).$$
We will also suppress some terms in the estimator $\widetilde{f}_{\ell}$ in order to get stochastic independence between 
the two couples of random functions $(m_{i\ell},\sigma_{i\ell})$ and $(m_{j\ell},\sigma_{j\ell})$ involved in the U-statistic.  
To this end we set for $\theta\in\Theta$, $v\in I$ and $1\leq i,j\leq n$, 
$$A_{v,ij}(\theta)=\frac{1}{\sigma_{i\ell}(\theta)}K_b\left[L_{v,i\ell}(\theta)-\ep_{j\ell}(\theta)\right].$$
Using the condition $\frac{\ell}{\sqrt{n}b}=o(1)$, we have
$$\sup_{v\in I}\left\vert \widetilde{f}_{\ell}(v)-\frac{1}{n^2}\sum_{(i,j)\in\mathcal{I}_n}A_{v,ij}(\hat{\theta})\right\vert=o_{\P}(1).$$

\subsubsection{Outline of the proof}
The goal of the proof is to show that 
\begin{equation}\label{superbut1}
\sup_{v\in I}\left|\frac{1}{n^{3/2}}\sum_{1\leq i,j\leq n}\left[A_{v,ij}(\hat{\theta})-A_{v,ij}(\theta_0)
-\dot{A}_{v,ij}(\theta_0)^T (\hat{\theta}-\theta_0)\right]\right|=o_{\P}(1)
\end{equation}
and in a second step that
\begin{equation}\label{superbut2}
\sup_{v\in I}\left\vert \frac{1}{n^2}\sum_{(i,j)\in\mathcal{I}_n}\dot{A}_{v,ij}(\theta_0)-\dot{h}_{\theta_0}(v)\right\vert=o_{\P}(1).
\end{equation}
In the proof of Theorem $1$, we have already shown that 
$$\sqrt{n}\sup_{v\in I}\left\vert \check{f}_X(v)-\frac{1}{n^2}\sum_{(i,j)\in\mathcal{I}_n}A_{v,ij}(\theta_0)\right\vert=o_{\P}(1).$$
Note that from assumption {\bf A4}, assertion (\ref{superbut1}) will hold if for all $M>0$ and integers $n$ such that $M/\sqrt{n}<\epsilon$, we have
$$\sup_{v\in I, \theta\in \Theta_{0,n}}\left|\frac{1}{n^2}\sum_{1\leq i,j\leq n}\left[A_{v,ij}(\theta)-A_{v,ij}(\theta_0)
-\dot{A}_{v,ij}(\theta_0)^T (\theta-\theta_0)\right]\right|=o_{\P}(1/\sqrt{n}),$$
where $\Theta_{0,n}$ is a short notation for $\Theta_{0,M/\sqrt{n}}$.
 We will show the following sufficient condition 
$$\sup_{v\in I, \theta\in \Theta_{0,n}}\left\|\frac{1}{n^2}\sum_{1\leq i,j\leq n}\left[\dot{A}_{v,ij}(\theta)-\dot{A}_{v,ij}(\theta_0)\right]\right\|=o_{\P}(1).$$
Then the two assertions (\ref{superbut1}) and (\ref{superbut2}) (and then Theorem $2$) will follow if we show that 
\begin{equation}\label{efficace1}
\frac{1}{n^2}\sup_{v\in I, \theta\in \Theta_{0,n}}\left\Vert\sum_{(i,j)\in\mathcal{I}_n}\left[\dot{A}_{v,ij}(\theta)-\E_{Y_{j\ell}}\dot{A}_{v,ij}(\theta)\right]\right\Vert=o_{\P}(1),
\end{equation}
\begin{equation}\label{efficace2}
\frac{1}{n^2}\sup_{v\in I, \theta\in \Theta_{0,n}}\left\Vert\sum_{(i,j)\in\mathcal{I}_n}\left[\E_{Y_{j\ell}}\dot{A}_{v,ij}(\theta)-\E_{Y_{j\ell}}\dot{A}_{v,ij}(\theta_0)\right]\right\Vert=o_{\P}(1)
\end{equation}
and 
\begin{equation}\label{efficace3}
\sup_{v\in I}\left\Vert \frac{1}{n^2}\sum_{(i,j)\in\mathcal{I}_n}\E_{Y_{j\ell}}\dot{A}_{v,ij}(\theta_0)-\dot{h}_{\theta_0}(v)\right\Vert=o_{\P}(1),
\end{equation}
where the function $h_{\theta}$ is defined before the statement of Theorem $2$.
Assertion (\ref{efficace1}) will be studied using martingale properties (see the subsection \ref{martii}).
In the rest of the proof, we prove the assertions (\ref{efficace1}), (\ref{efficace2}) and (\ref{efficace3}).  
Note that we have the following expression.
$$\dot{A}_{v,ij}(\theta)=\dot{\sigma^{-1}}_{i\ell}(\theta)K_b\left[L_{v,i\ell}(\theta)-\ep_{j\ell}(\theta)\right]+\frac{\dot {L}_{v,i\ell}(\theta)-\dot{\ep}_{j\ell}(\theta)}{\sigma_{i\ell}(\theta)}K_b'\left[L_{v,i\ell}(\theta)-\ep_{j\ell}(\theta)\right].$$

\subsubsection{Proof of assertion (\ref{efficace1})}
We set 
$$S_n(v,\theta)=\frac{1}{n^2}\sum_{(i,j)\in\mathcal{I}_n}\left[\dot{A}_{v,ij}(\theta)-\E_{Y_{j\ell}}\left(\dot{A}_{v,ij}(\theta)\right)\right].$$
Let $\eta=\eta_n$ a sequence of positive real numbers such that $\eta/b^3=o_{\P}(1)$. We take for instance $\eta=n^{-\frac{3}{3+\delta}}$. 
Let $\left\{(v_h,\theta_h):h\in \mathcal{H}\right\}$ a family of points in $I\times \Theta_{0,n}$ such that 
for $(v,\theta)\in I\times \Theta_{0,n}$, there exists $h\in\mathcal{H}$ such that $\max\left\{\vert v-v_h\vert,\Vert\theta-\theta_h\Vert\right\}\leq \eta$. The set $\mathcal{H}$ can be chosen such that $\left\vert \mathcal{H}\right\vert=O\left(\eta^{-d-1}\right)=O\left(n^{\frac{3(d+1)}{3+\delta}}\right)$. 
Using Lemma \ref{propres2} (3), we first notice that
$$\sup_{(v,\theta)\in I\times \Theta_{0,n}}\left\Vert S_n(v,\theta)\right\Vert-\sup_{(v,\theta)\in \mathcal{H}}\left\Vert S_n(v,\theta)\right\Vert=O_{\P}\left(\frac{\eta}{b^3}\right)=o_{\P}(1).$$
To show (\ref{efficace1}), it remains to prove that 
$$\sup_{(v,\theta)\in \mathcal{H}}\left\Vert S_n(v,\theta)\right\Vert=o_{\P}(1).$$
But this is a consequence of Lemma \ref{martin} applied coordinatewise to $S_n(\cdot,\cdot)$, using the martingale decomposition (\ref{decompmart}). The assumptions used in Lemma  
\ref{martin} can checked using Lemma \ref{propres2}.

\subsubsection{Proof of assertion (\ref{efficace2})}
The following notations will be needed. We define 
$$\Gamma_i^{(1)}(v,\theta)=\dot{\sigma_{i\ell}^{-1}}(\theta)\int K(w)f_{\theta}\left(L_{v,i\ell}(\theta)-bw\right)dw,$$
$$\Gamma_i^{(2)}(v,\theta)=\frac{\dot{L}_{v,i\ell}(\theta)}{\sigma_{i\ell}(\theta)}\int K(w)f'_{\theta}\left(L_{v,i\ell}(\theta)-bw\right)dw,$$
$$\Gamma_i^{(3)}(v,\theta)=\frac{1}{\sigma_{i\ell}(\theta)}\int K(w)\dot{f}_{\theta}\left(L_{v,i\ell}(\theta)-bw\right)dw.$$
Note that from Lemma \ref{propres}, $4.$, we have 
$$\E_{Y_{j\ell}}\dot{A}_{v,ij}(\theta)=\sum_{h=1}^3\Gamma_i^{(h)}(v,\theta).$$
Then assertion (\ref{efficace2}) will follow if we show that for $h=1,2,3$, 
\begin{equation}\label{premiere}
\frac{1}{n^2}\sum_{(i,j)\in\mathcal{I}_n}\sup_{v\in I,\theta\in \Theta_{0,n}}\left\|\Gamma_i^{(h)}(v,\theta)-\Gamma_i^{(h)}(v,\theta_0)\right\|=o_{\P}(1).
\end{equation}
The proof of (\ref{premiere}) follows from the following bounds, Assumption {\bf A3}, Lemma \ref{propres} and Lemma \ref{propres3}. 
\begin{itemize}
\item For $h=1$, we have
\begin{eqnarray*}
&&\sup_{v\in I,\theta\in\Theta_{0,n}}\left\|\Gamma_i^{(1)}(v,\theta)-\Gamma_i^{(1)}(v,\theta_0)\right\|\\
&\leq& \frac{C}{\sqrt{n}}\left[\nor{\ddot{\sigma^{-1}_{i\ell}}}+\nor{\dot{\sigma^{-1}_{i\ell}}}\left(\sup_{v\in I}\nor{\dot{L}_{v,i\ell}}+C_1+C_2\left(\sup_{v\in I}\nor{L_{v,i\ell}}+1\right)\right)\right],
\end{eqnarray*}
where $C_1$ and $C_2$ are the constants given in Lemma \ref{propres} (3).
\item
For $h=2$, we have 
\begin{eqnarray*}
\sup_{v\in I,\theta\in\Theta_{0,n}}\left\|\Gamma_i^{(2)}(v,\theta)-\Gamma_i^{(2)}(v,\theta_0)\right\|
&\leq& \frac{C}{\sqrt{n}}\left[\sup_{v\in I}\nor{\ddot{L}_{v,i\ell}}+\sup_{v\in I}\nor{\dot{L}_{v,i\ell}}\cdot\nor{\dot{\sigma}_{i\ell},\sigma_{i\ell}}\right]\\
&+&\frac{C}{\sqrt{n}}\sup_{v\in I}\nor{\dot{L}_{v,i\ell}}\cdot\left(1+\sup_{v\in I}\nor{L_{v,i\ell}}+\sup_{v\in I}\nor{\dot{L}_{v,i\ell}}\right).
\end{eqnarray*}

\item
Finally we have
\begin{eqnarray*}
\sup_{v\in I,\theta\in\Theta_{0,n}}\left\|\Gamma_i^{(3)}(v,\theta)-\Gamma_i^{(3)}(v,\theta_0)\right\|
&\leq& \frac{C}{\sqrt{n}}\left[\nor{\dot{\sigma}_{i\ell},\sigma_{i\ell}}\left(1+\sup_{v\in I}\nor{L_{v,i\ell}}\right)\right]\\
&+&\frac{C}{\sqrt{n}}\left[1+\sup_{v\in I}\nor{L_{v,i\ell}}^2+\sup_{v\in I}\nor{\dot{L}_{v,i\ell}}\cdot\left(1+\sup_{v\in I}\nor{L_{v,i\ell}}\right)\right].
\end{eqnarray*}
Using the integrability properties stated in Lemma \ref{propres3} and assumption {\bf A3}, assertion (\ref{efficace2}) follows.
\end{itemize}
\subsubsection{Proof of assertion (\ref{efficace3})}
We set
$$\Delta_{v,i}^{(1)}=\frac{-\dot{\sigma}_{i\ell}(\theta_0)}{\sigma_{i\ell}^2(\theta_0)}f_{\theta_0}\left(L_{v,i\ell}(\theta_0)\right),
\quad \Delta_{v,i}^{(2)}=\frac{\dot{L}_{v,i\ell}(\theta_0)}{\sigma_{i\ell}(\theta_0)}f'_{\theta_0}\left(L_{v,i\ell}(\theta_0)\right),
\quad \Delta_{v,i}^{(3)}=\frac{1}{\sigma_{i\ell}(\theta_0)}\dot{f}_{\theta_0}\left(L_{v,i\ell}(\theta_0)\right).$$
Using Lemma \ref{propres} and Lemma \ref{propres3}, it is easily seen that for $h=1,2,3$,
$$\frac{1}{n^2}\sum_{(i,j)\in\mathcal{I}_n}\sup_{v\in I}\left\Vert \Gamma_{v,i}(\theta_0)-\Delta_{v,i}\right\Vert=o_{\P}(1).$$
To end the proof of assertion (\ref{efficace3}), it remains to show that for $h=1,2,3$, 
\begin{equation}\label{deuxieme}
\sup_{v\in I}\frac{1}{n^2}\left\vert\sum_{(i,j)\in\mathcal{I}_n}\left[\Delta_{v,i}^{(h)}-\E\left(\Delta_{v,i}^{(h)}\right)\right]\right\vert =o_{\P}(1).
\end{equation}
Note that $\sum_{h=1}^3\E\left(\Delta_{v,i}^{(h)}\right)=\dot{h}_{\theta_0}(v)$.
To show (\ref{deuxieme}), we first notice that $\sup_{n\in \N^{*},v\in I}\E\Vert \Delta_{v,i}^{(h)}\Vert<\infty$ for $h=1,2,3$. 
Then, since $\ell=o(n)$, it is easily seen that 
$$\sup_{v\in I}\left\vert\frac{1}{n^2}\sum_{(i,j)\in\mathcal{I}_n}\left[\Delta_{v,i}^{(h)}-\E\left(\Delta_{v,i}^{(h)}\right)\right]
-\frac{1}{n}\sum_{i=1}^{n'}\left[\Delta_{v,i}^{(h)}-\E\left(\Delta_{v,i}^{(h)}\right)\right]\right\vert=o_{\P}(1).$$
Now, we set for $v\in I$ and $h=1,2,3$, 
$$G^{(h)}_n(v)=\frac{1}{n}\sum_{i=1}^{n'}\left[\Delta_{v,i}^{(h)}-\E\left(\Delta_{v,i}^{(h)}\right)\right].$$ 
Then assertion (\ref{efficace3}) will follow if we show that
\begin{equation}\label{onva}
\sup_{v\in I}\left\Vert G^{(h)}_n(v)\right\Vert=o_{\P}(1),\quad h=1,2,3.
\end{equation}
We have using the $\ell-$dependence,
$$\E\left\vert G^{(h)}_n(v)\right\vert^2\leq \frac{\ell}{n^2}\sum_{s=1}^{\ell}\sum_{g=0}^{k-1}\v\left(M^{(h)}_{s+g\ell}{(h)}\right).$$
Moreover, using the fact that $\ell=o(n)$ and that $\sup_{n,i\geq 1}\sup_{v\in I} \v\left(M^{(h)}_i(v)\right)$ is bounded, we get  $G_n^{(h)}(v)=o_{\P}(1)$ for each $v\in I$.
Now, (\ref{onva}) will follow if we show that $\sup_{v\neq \bar{v}}\frac{\Vert G_n^{(h)}(v)-G_n^{(h)}(\bar{v})\Vert}{\vert v-\bar{v}\vert}=O_{\P}(1)$. But this is a consequence of the following bounds. First, there exists $C>0$ such that
$$\left\Vert\Delta^{(1)}_{v,i}-\Delta^{(1)}_{v',i}\right\Vert\leq C \nor{\dot{\sigma}_{i\ell},\sigma_{i\ell}}\vert v-v'\vert,$$
$$\left\Vert\Delta^{(2)}_{v,i}-\Delta^{(2)}_{v',i}\right\Vert\leq C\left[\nor{\dot{\sigma}_{i\ell},\sigma_{i\ell}}+\sup_{v\in I}\nor{\dot{L}_{v,i\ell}}\right]\cdot\vert v-v'\vert.$$
For $h=3$, we set $E_1=\E\left[\frac{\dot{\sigma}_j(\theta_0)}{\sigma_j(\theta_0)}\right]$ and 
$E_2=\E\left[\frac{\dot{m}_j(\theta_0)}{\sigma_j(\theta_0)}\right]$. Then $\dot{f}_{\theta_0}(w)=E_1\left(f_{\ep}(w)+wf'_{\ep}(w)\right)+E_2 f'_{\ep}(w)$. Then, using assumption {\bf A5}, we have 
\begin{eqnarray*}
&&\left\Vert\Delta^{(3)}_{v,i}-\Delta^{(3)}_{v',i}\right\Vert\\
&\leq& C\left[\vert v-v'\vert+\left\vert L_{v,i\ell}(\theta_0)f'_{\ep}\left(L_{v,i\ell}(\theta_0)\right)-L_{v',i\ell}(\theta_0)f'_{\ep}\left(L_{v',i\ell}(\theta_0)\right)\right\vert\right]\\
&\leq& C\vert v-v'\vert\cdot\left[1+\sup_{v\in I}\nor{L_{v,i\ell}}\right]. 
\end{eqnarray*}
In the previous bounds, the real number $C$ does not depends on $v,\bar{v}\in I$. Then (\ref{onva}) follows and the proof of assertion (\ref{efficace3}) is now complete.

 This ends the proof of Theorem $2$.$\square$

\subsection{Proof of Corollary $1$}
From Theorem $1$, Theorem $2$ and assumption {\bf A8}, we have 
$$\sqrt{n}\left[\hat{f}_X(v)-f_X(v)\right]=\frac{1}{\sqrt{n}}\sum_{i=1}^n \overline{M}_{i,v}+o_{\P}(1).$$
The first part of the corollary concerns the convergence of the finite dimensional distributions. Note also that the previous convergence is uniform if assumption {\bf A7} holds true. Convergence of finite dimensional distributions is straightforward using a central lime theorem for weakly dependent time series. For instance, the central limit theorem given in \citet{Zh}, Theorem $3$, applies in our case. For the uniform convergence, it remains to show the tightness of the empirical process $G_n:v\mapsto \frac{1}{\sqrt{n}}\sum_{i=1}^n\overline{M}_{i,v}$.
Since we assumed {\bf A7}, it is only necessary to study the tightness in $\mathcal{C}(I)$, the space of real-valued and continuous function defined on $I$, of 
$$\widetilde{G}_n:v\mapsto \frac{1}{\sqrt{n}}\sum_{i=1}^n \left[\frac{1}{\sigma_i}f_{\ep}\left(\frac{v-m_i}{\sigma_i}\right)-\E \frac{1}{\sigma_i}f_{\ep}\left(\frac{v-m_i}{\sigma_i}\right)\right].$$
To this end, we use the Kolmogorov-Chentsov criterion (see for instance \citet{Kal}, Corollary $14.9$). If $v,\bar{v}\in I$ satisfy $v\leq \bar{v}$, we have from Jensen inequality, 
\begin{eqnarray*}
\E\left\vert \widetilde{G}_n(v)-\widetilde{G}_n(v)\right\vert^2&\leq& \vert v-\bar{v}\vert \int \E\left\vert \widetilde{G}'_n(u)\right\vert^2 du\\
&\leq& \vert v-\bar{v}\vert^2\sup_{u\in I}\E\left\vert \widetilde{G}'_n(u)\right\vert^2.
\end{eqnarray*}
Then the tightness will follow if we show that 
\begin{equation}\label{fin}
\sup_{u\in I}\E\left\vert \widetilde{G}'_n(u)\right\vert^2=O(1).
\end{equation}
Note that 
$$\widetilde{G}'_n(u)=\frac{1}{\sqrt{n}}\sum_{i=1}^n \left[\frac{1}{\sigma_i^2}f_{\ep}'\left(\frac{v-m_i}{\sigma_i}\right)
-\E\frac{1}{\sigma_i^2}f_{\ep}'\left(\frac{v-m_i}{\sigma_i}\right)\right].$$
Moreover, for $i\leq j$, $v\in I$ and $\ell=j-i$, we have
\begin{eqnarray*}
&&\cov\left[\frac{1}{\sigma_i^2}f_{\ep}'\left(\frac{v-m_i}{\sigma_i}\right),\frac{1}{\sigma_j^2}f_{\ep}'\left(\frac{v-m_j}{\sigma_j}\right)\right]\\
&\leq& \frac{\Vert f_{\ep}'\Vert_{\infty}}{\gamma^2}\E\left\vert \frac{1}{\sigma_j^2}f_{\ep}'\left(\frac{v-m_j}{\sigma_j}\right)-\frac{1}{\sigma_{j\ell}^2}f_{\ep}'\left(\frac{v-m_{j\ell}}{\sigma_{j\ell}}\right)\right\vert.
\end{eqnarray*}
Then using assumption {\bf A2} and assumption {\bf A5}, we deduce that there exist $C>0$ such that for all $i\leq j$, 
$$\cov\left[\frac{1}{\sigma_i^2}f_{\ep}'\left(\frac{v-m_i}{\sigma_i}\right),\frac{1}{\sigma_j^2}f_{\ep}'\left(\frac{v-m_j}{\sigma_j}\right)\right]\leq C a^{\frac{j-i}{2}},$$
where $a\in (0,1)$ is defined in assumption {\bf A2}. The proof of the last inequality uses the same arguments than the proof of Lemma \ref{propres3}. This control of covariances immediately implies (\ref{fin}). The tightness criterion of Kolmogorov-Chentsov applies. The proof of Corollary $1$ is now complete.$\square$

\subsection{Checking the regularity assumptions on densities}

\subsubsection*{Density regularities of ARCH processes}

Here, we assume that $(X_t)$ is a stationnary ARCH process defined by 
$$X_t=\ep_t\sigma_t,\quad \sigma^2_t=\alpha_0+\sum_{j\geq 1}\alpha_jX^2_{t-j}.$$
We assume here that $\sum_{j=1}^{\infty}\alpha_j<\infty$ and that $f_{\ep}$ is bounded. We set $\mu(dx)=\sup_{\vert z\vert\leq z_0}f_{\ep}(x+z)dx$.  
We denote by $f_{\sigma^2}$ the probability density of the conditional variance $\sigma_t^2$ and for $v\neq 0$, $g_v(x)=\frac{2v}{x^2}f_{\sigma^2}\left(\frac{v^2}{x^2}\right)$ which is well defined for $x\neq 0$.
We also set $G(x)=\sup_{v\in I}g_v(x)$ and for a compact interval $I$ which does not contain $0$, $\mathcal{G}_I=\left\{g_v: v\in I\right\}$.
\begin{lem}\label{regule}
\begin{enumerate}
\item
Assume that $\alpha_1,\alpha_2>0$. Then $f_{\sigma^2}$ is bounded. Moreover, if $\E\sigma_t<\infty$, then for all $v\neq 0$,
we have $\int g_v(x)^2\mu(dx)<\infty$.
\item
Assume that $\alpha_1,\alpha_2,\alpha_3>0$. Then there exists a constant $C>0$ such that for $s_1,s_2\geq 0$, we have 
$$\left\vert f_{\sigma^2}(s_2)-f_{\sigma^2}(s_1)\right\vert\leq C \left(1+\sqrt{\vert s_1\vert}\wedge \sqrt{\vert s_2\vert}\right)\cdot \sqrt{\vert s_2-s_1\vert}.$$ 
\item
In addition to the previous point, assume that there exists $\delta\in (0,\frac{1}{2})$ such that $\E\sigma_t^{\frac{3}{2}+\delta}<\infty$ and that $x\mapsto \vert x\vert^{\frac{3}{2}+\delta} f_{\ep}(x)$ is bounded. Then
there exists a real number $o>0$ such that $\int G(x)^{2+o}\mu(dx)<\infty$. Moreover there exists some constants $\zeta,C>0$ such that 
$$N_{[]}\left(\epsilon,\mathcal{G}_I,\L^2(\mu)\right)\leq C \epsilon^{-\zeta}.$$
\end{enumerate}
\end{lem}

\paragraph{Proof of Lemma \ref{regule}.}
Before proving the lemma, we first derive an expression for $f_{\sigma^2}$ involving conditional distributions. 
We will use for $j\geq 1$ the notation ${\bf z_j}=(z_j,z_{j+1},\ldots)$, we set 
We set $k_3({\bf z_3})=\alpha_0+\sum_{j=3}^{\infty}\alpha_j z_j^2$ and 
$$s({\bf z_1})=\sqrt{\alpha_0+\sum_{j=1}^{\infty}\alpha_j z_j^2}.$$
The measure $\nu$ will denote the probability distribution of $\left(X_{t-1},X_{t-2},\ldots\right)$. 
Moreover, we set 
$$r(z_1,z_2\vert {\bf z_3})=\frac{1}{s({\bf z_2})}f_{\ep}\left(\frac{z_1}{s({\bf z_2})}\right)\frac{1}{s({\bf z_3})}f_{\ep}\left(\frac{z_2}{s({\bf z_3})}\right)$$
and
$$\bar{r}(z_1,z_2\vert {\bf z_3})=\frac{1}{\sqrt{\alpha_1\cdot\alpha_2}}r\left(\frac{z_1}{\sqrt{\alpha_1}},\frac{z_2}{\sqrt{\alpha_2}}\big\vert {\bf z_3}\right).$$
If $h:\R\rightarrow\R$ is a bounded and measurable function, we have
\begin{eqnarray*}
\E h\left(\sigma_t^2\right)&=& \int\int\int h\left(\alpha_1z_1^2+\alpha_2z_2^2+k_3({\bf z_3})\right)r(z_1,z_2\vert {\bf z_3})dz_1dz_2d\nu({\bf z_3})\\
&=& \int\int\int h\left(z_1^2+z_2^2+k_3({\bf z_3})\right)\bar{r}(z_1,z_2\vert {\bf z_3})dz_1dz_2d\nu({\bf z_3})\\
&=& \int_{0}^{\infty}\int_{-\pi}^{\pi}\int h\left(\rho^2+k_3({\bf z_3})\right)\bar{r}(\rho\cos(\phi),\rho\sin(\phi)\vert {\bf z_3})\rho d\rho d\phi d\nu({\bf z_3})\\
&=& \int \int_{k_3({\bf z_3})\leq s}\int_{-\pi}^{\pi} h(s)\frac{1}{2}\bar{r}(\sqrt{s-k_3({\bf z_3})}\cos(\phi),\sqrt{s-k_3({\bf z_3})}\sin(\phi)\vert {\bf z_3})ds d\phi d\nu({\bf z_3}).\\
\end{eqnarray*}
Then we deduce that for $s\geq 0$,
\begin{equation}\label{express}
f_{\sigma^2}(s)=\frac{1}{2}\int_{k_3({\bf z_3})\leq s}\int_{-\pi}^{\pi}\bar{r}(\sqrt{s-k_3({\bf z_3})}\cos(\phi),\sqrt{s-k_3({\bf z_3})}\sin(\phi)\vert {\bf z_3})d\phi d\nu({\bf z_3}).
\end{equation}

\begin{enumerate}
\item
The fact that $f_{\sigma^2}$ is bounded is a consequence of the expression (\ref{express}), using the fact that $f_{\ep}$
and then $\bar{r}$ are bounded.
Moreover, we have 
\begin{eqnarray*}
\int g_v(x)^2\mu(dx)&\leq& 4v^2\cdot\Vert f_{\ep}\Vert_{\infty}\int \frac{1}{x^4}f_{\sigma^2}\left(\frac{v^2}{x^2}\right)dx\\
&\leq& \frac{4}{\vert v\vert}4\Vert f_{\ep}\Vert_{\infty}\int_0^{\infty}\sqrt{y}f_{\sigma^2}(y)dy.
\end{eqnarray*}
This bound gives the result.
\item
We use the expression (\ref{express}). Using some basic computations, it is easily seen that for real numbers $z_1,z_2$,
$$\left\vert r(z_1,z_2\vert {\bf z_3})-r(\bar{z}_1,\bar{z}_2\vert {\bf z_3})\right\vert\leq C\left(\vert z_1-\bar{z}_1\vert+\left(1+\vert z_1\vert\wedge \vert \bar{z}_1\vert\right)\cdot\vert z_2-\bar{z}_2\vert\right)$$
for some constant $C>0$. Setting now for $s\geq k_3({\bf z_3})$,
$$h(s,\phi,{\bf z_3})=\bar{r}\left[\sqrt{s-k_3({\bf z_3})}\cos(\phi),\sqrt{s-k_3({\bf z_3})}\sin(\phi)\vert {\bf z_3}\right],$$
we have for $s_2\geq s_1\geq k_3({\bf z_3})$,
$$\left\vert h(s_1,\phi,{\bf z_3})-h(s_2,\phi,{\bf z_3})\right\vert\leq C\left(1+\sqrt{s_1\wedge s_2}\right)\cdot\sqrt{\vert s_2-s_1\vert}.$$
We deduce that for $s_1,s_2\geq 0$,
$$\left\vert f_{\sigma^2}(s_2)-f_{\sigma^2}(s_1)\right\vert\leq  C\left(1+\sqrt{s_1\wedge s_2}\right)\cdot\sqrt{\vert s_2-s_1\vert}+\P\left(s_1<k_3(X_{t-3},\ldots)\leq s_2\right).$$

Moreover, it is easily seen that $k_3\left(X_{t-3},\ldots\right)$ has a density $q_k$ such that 
$$q_k(x)\leq C\int_{k_4({\bf z_4})\leq x}\left(x-k_4({\bf z_4})\right)^{-1/2}d\nu({\bf z_4}),$$
where $k_4({\bf z_4})=\alpha_0+\sum_{j\geq 4}\alpha_jz_j^4$.
Then, it can be shown that 
$$\P\left(s_1<k_3(X_{t-3},\ldots)\leq s_2\right)\leq C\sqrt{\vert s_2-s_1\vert}.$$
This proves the second point of this lemma.

\item
We first show that under our assumptions,
\begin{equation}\label{decreas}
\sup_{x>0}x^{\frac{3}{4}+\frac{\delta}{2}}f_{\sigma^2}(x)<\infty. 
\end{equation}
We first observe that 
\begin{eqnarray*}
&&\vert u\vert^{\frac{3}{2}+\delta}\bar{r}\left(u\cos\phi,u\sin\phi\vert{\bf z_3}\right)\\
&\leq& C \left[\vert u\cos\phi\vert^{\frac{3}{2}+\delta}+\vert u\sin\phi\vert^{\frac{3}{2}+\delta}\right]\cdot \bar{r}\left(u\cos\phi,u\sin\phi\vert{\bf z_3}\right)\\
&\leq& C\sup_{y>0}\left\{y^{\frac{3}{2}+\delta}f_{\ep}(y)\right\}\cdot\frac{\left\vert\alpha_0+\alpha_1u^2\sin^2\phi+\sum_{j\geq 2}\alpha_jz_{j+1}^2\right\vert^{\frac{1}{4}+\frac{\delta}{2}}}{s({\bf z_3})^{\frac{1}{2}+\delta}}f_{\ep}\left(\frac{u\sin\phi}{\sqrt{\alpha_2} s({\bf z_3})}\right)\\
&+& C \sup_{y>0}\left\{y^{\frac{3}{2}+\delta}f_{\ep}(y)\right\}\cdot s({\bf z_3})^{\frac{1}{2}+\delta}\\
&\leq& C\sup_{y>0}\left\{y^{\frac{3}{2}+\delta}f_{\ep}(y)\right\}\cdot\left[\sup_{y>0}\left\{y^{\frac{1}{2}+\delta}f_{\ep}(y)\right\}+\left\vert\alpha_0+\sum_{j\geq 2}\alpha_j z_{j+1}^2\right\vert^{\frac{1}{4}+\frac{\delta}{2}}+s({\bf z_3})^{\frac{1}{2}+\delta}\right]\\
\end{eqnarray*}
Then we deduce that the function $\vartheta$ defined by 
$$\vartheta({\bf z_3})=\sup_{u\in\R,\phi\in(-\pi,\pi)}\vert u\vert^{\frac{3}{2}+\delta}\bar{r}\left(u\cos\phi,u\sin\phi\vert{\bf z_3}\right)$$
satisfies $\E\left[\vartheta\left(X_{t-3},X_{t-4},\ldots\right)\right]<\infty$.
Using the expression (\ref{express}), condition $\E\sigma_t^{\frac{3}{2}+\delta}<\infty$ and the decomposition $x=x-k_3({\bf z_3})+k_3({\bf z_3})$, we get the bound 
\begin{eqnarray*}
&&x^{\frac{3}{4}+\frac{\delta}{2}}f_{\sigma^2}(x)\\
&\leq& C\left(\E\left[k_3(X_{t-3},X_{t-4},\ldots)^{\frac{3}{4}+\frac{\delta}{2}}\right]+\E\left[\vartheta\left(X_{t-3},X_{t-4},\ldots\right)\right]\right)\\
&<&\infty.
\end{eqnarray*}
This shows (\ref{decreas}).
For simplicity, we now assume that $I\subset (0,\infty)$ (the case $I\subset (-\infty,0)$ is identical). 
First, we show that there exists $o>0$ such that $\int G(x)^{2+o}\mu(dx)<\infty$. We have 
Since $I$ is compact and $\sigma^2_t(\theta_0)$ is bounded from below, there exists $\varpi>0$ such that 
$G(x)=0$ when $x>\varpi$. We choose $o$ such that $\frac{3}{4}+\frac{\delta}{2}=1-\frac{1-o}{2(2+o)}$. We get 
\begin{eqnarray*}
\int_0^{\varpi}G(x)^{2+o}dx&\leq& C\int_0^{\varpi}\frac{1}{x^{1-o}}\left\vert \frac{x^{\frac{1-o}{2+o}}}{x^2}\sup_{v\in I}f_{\sigma^2}\left(\frac{v^2}{x^2}\right)\right\vert^{2+o}dx\\
&\leq& C\left[\sup_{y>0}\left\{y^{\frac{3}{4}+\frac{\delta}{2}}f_{\sigma^2}(y)\right\}\right]^{2+o}.
\end{eqnarray*} 
This shows that $\int_0^{\infty}G(x)^{2+o}dx<\infty$. Finally, we consider the bracketing numbers of the family $\mathcal{G}_I$. Let $\eta\in (0,1)$ be such that 
$8\eta +\left(\frac{1}{2}-\delta\right)\left(2-2\eta\right)<1$. 
Let $v_1,v_2\in I$.
We have 
\begin{eqnarray*}
\left\vert g_{v_1}(x)-g_{v_2}(x)\right\vert&\leq& C G(x)\vert v_1-v_2\vert+\frac{2v}{x^2}\left\vert f_{\sigma^2}\left(\frac{v_1^2}{x^2}\right)-f_{\sigma^2}\left(\frac{v_2^2}{x^2}\right)\right\vert\\
&\leq& C\left[G(x)\vert v_1-v_2\vert+\frac{G(x)^{1-\eta}}{x^{2\eta}}\left\vert f_{\sigma^2}\left(\frac{v_1^2}{x^2}\right)-f_{\sigma^2}\left(\frac{v_2^2}{x^2}\right)\right\vert^{\eta}\right]\\
&\leq& C\left[G(x)+G(x)^{1-\eta}\left(\vert x\vert^{-3\eta}+x^{-4\eta}\right)\right]\cdot\vert v_1-v_2\vert^{\eta}.
\end{eqnarray*}
Since we have
$$\int_0^{\varpi}G(x)^{2+o}dx<\infty,\quad \int_0^{\varpi}G(x)^{2-2\eta}x^{-8\eta}dx<\infty,$$
where the second integrability condition follows from the assumption on $\eta$ and (\ref{decreas}), the bound given 
for $N_{[]}\left(\epsilon,\mathcal{G}_I,\L^2(\mu)\right)$ easily follows (see \cite{VW}, Example $19.7$, for $\eta=1$ and $\mu$ a probability measure, but the arguments are similar in our case). This completes the proof of the lemma.$\square$

\end{enumerate}

\subsubsection*{A result for ARMA-GARCH processes}

\begin{lem}\label{ARMAGARCH}
Assume that $(\psi_j)_{j\geq 1}$ and $(\alpha_j)_{j\geq 0}$ are two summable sequences of real numbers such that $\alpha_i\geq 0$ for $i\geq 0$ and $\alpha_i>0$ for $0\leq i\leq 3$ and there exists an integer $q\geq 1$ such that $\psi_q\neq 0$. Let $\left(Z_t\right)_{t\in\Z}$ be a stationary process of real random variables such that $\E Z_t^2<\infty$ and such that the conditional distribution of $Z_t\vert Z_{t-1},Z_{t-2},\ldots$ has a bounded density. 
Then the density $\omega$ of the couple $\left(\sum_{j=1}^{\infty}\psi_j Z_{t-j},\sqrt{\alpha_0+\sum_{j=1}^{\infty}\alpha_j Z_{t-j}^2}\right)$ satisfies $\omega(x,y)\leq C y$ for a positive constant $C$.
\end{lem}

\paragraph{Proof of Lemma \ref{ARMAGARCH}.}
We set ${\bf z_i}=\left(z_i,z_{i+1},\ldots\right)$ and we denote by $\widetilde{f}\left(\cdot\vert {\bf z_1}\right)$ the conditional density of $Z_t$ given that 
$Z_{t-i}=z_i$ for $i\geq 1$.
We consider two cases.
\begin{enumerate}
\item
We first assume that $q\leq 3$. We set 
$$k_1({\bf z_4})=\sum_{j\geq 4}\psi_j z_j,\quad k_2({\bf z_4})=\alpha_0+\sum_{j\geq 4}\alpha_j z_j^2$$
and $r=\Vert a\Vert$ with $a=\left(\frac{\psi_i}{\sqrt{\alpha_i}}\right)_{1\leq i\leq 3}$.
We also set 
$$\zeta({\bf z_1})=\frac{1}{\sqrt{\alpha_1\alpha_2\alpha_3}}\widetilde{f}\left(\frac{z_1}{\sqrt{\alpha_1}}\big\vert \frac{z_2}{\sqrt{\alpha_2}},\frac{z_3}{\sqrt{\alpha_3}},{\bf z_4}\right)\cdot\widetilde{f}\left(\frac{z_2}{\sqrt{\alpha_2}}\big\vert \frac{z_3}{\sqrt{\alpha_3}},{\bf z_3}\right)\cdot\widetilde{f}\left(\frac{z_3}{\sqrt{\alpha_3}}\big\vert{\bf z_4}\right).$$ 
We also denote by $\mathcal{R}$ the rotation of $\R^3$ such that $\mathcal{R}e_1=a/r$ where $e_1=(1,0,0)$. 
For simplicity of notations, we only use one sign "integral" and do not precise the boundaries for integration in the next computations.
Then if $\nu$ denotes the probability distribution of $(Z_{t-4},Z_{t-5},\ldots)$, we have 
\begin{eqnarray*}
&&\E h\left(\sum_{j=1}^{\infty}\psi_j Z_{t-j},\sqrt{\alpha_0+\sum_{j=1}^{\infty}\alpha_j Z_{t-j}^2}\right)\\
&=& \int h\left(\sum_{j=1}^3 a_j z_j+k_1({\bf z_4}),\sqrt{\sum_{j=1}^3 z_j^2+k_2({\bf z_4})}\right)\zeta({\bf z_1})dz_1dz_2dz_3d\nu({\bf z_4})\\
&=& \int h\left(r z_1+k_1({\bf z_4}),\sqrt{\sum_{j=1}^3z_j^2+k_2({\bf z_4})}\right)\zeta\left(\mathcal{R}(z_1,z_2,z_3),{\bf z_4}\right)dz_1dz_2dz_2d\nu({\bf z_4})\\
&=& \int h\left(r z_1+k_1({\bf z_4}),\sqrt{z_1^2+\rho^2+k_2({\bf z_4})}\right)\zeta\left(\mathcal{R}(z_1,\rho\cos\phi,\rho\sin\phi),{\bf z_4}\right)\rho dz_1d\rho d\phi d\nu({\bf z_4})\\
&=& \int h(x,y)\\&&\zeta\left(\mathcal{R}\left(\frac{x-k_1({\bf z_4})}{\alpha},\sqrt{y^2-\left[\frac{x-k_1({\bf z_4})}{\alpha}\right]-k_2({\bf z_4})}\cos\phi,\sqrt{y^2-\left[\frac{x-k_1({\bf z_4})}{\alpha}\right]-k_2({\bf z_4})}\sin\phi\right),{\bf z_4}\right)\\&&\frac{y}{r} dx dy d\phi d\nu({\bf z_4}).\\
\end{eqnarray*}
Since $\zeta$ is bounded, it is easily seen that $\omega(x,y)\leq C y$ for a positive constant $C$.
\item
We next consider the case $q\geq 4$. We set 
$$\zeta_1({\bf z_1})=\prod_{i=1}^2\widetilde{f}\left(z_i\vert {\bf z_{i+1}}\right),\quad \zeta_2({\bf z_3})=\prod_{i=3}^{q-1}\widetilde{f}\left(z_i\vert {\bf z_{i+1}}\right)$$
and $k_3({\bf z_3})=\alpha_0+\sum_{j=3}^{\infty}\alpha_j z_j^2$. 
For simplicity of notations, we assume that $\alpha_1=\alpha_2=1$  (otherwise, as in the previous point, a change of variables is needed in the computations given below).
Denoting by $\nu$, the probability distribution of $\left(Z_{q+1},Z_{q+2},\ldots\right)$, we have 
\begin{eqnarray*}
&&\E h\left(\sum_{j=1}^{\infty}\psi_j Z_{t-j},\sqrt{\alpha_0+\sum_{j=1}^{\infty}\alpha_j Z_{t-j}^2}\right)\\
&=& \int h\left(\sum_{j=q}^{\infty}\psi_j z_j,\sqrt{z_1^2+z_2^2+k_3({\bf z_3})}\right)\prod_{i=1}^q \widetilde{f}(z_i\vert {\bf z_{i+1}})dz_1\cdots dz_qd\nu({\bf z_{q+1}})\\
&=& \int h\left(\sum_{j=q}^{\infty}\psi_j z_j,\sqrt{\rho^2+k_3({\bf z_3})}\right)\zeta_1(\rho\cos\phi,\rho\sin\phi,{\bf z_3})\zeta_2({\bf z_3})\\
&&\widetilde{f}(z_q\vert {\bf z_{q+1}})\rho d\rho d\phi dz_3\cdots dz_q d\nu({\bf z_{q+1}})\\
&=& \int h\left(\sum_{j=q}^{\infty}\psi_j z_j,y\right)\zeta_1\left(\sqrt{y^2-k_3({\bf z_3})}\cos\phi, \sqrt{y^2-k_3({\bf z_3})}\cos\phi,{\bf z_3}\right)\\
&&\zeta_2({\bf z_3})\widetilde{f}(z_q\vert {\bf z_{q+1}})y dy d\phi dz_3\cdots dz_q d\nu({\bf z_{q+1}})\\
&=&\int h(x,y)\zeta_1\left(\sqrt{y^2-k_3({\bf z_3})}\cos\phi, \sqrt{y^2-k_3({\bf z_3})}\cos\phi,z_3,\ldots,z_{q-1},\frac{x-\sum_{j\geq q+1}\psi_j z_j}{\psi_q},{\bf z_{q+1}}\right)\\
&&\zeta_2\left(z_3,\ldots,z_{q-1},\frac{x-\sum_{j\geq q+1}\psi_j z_j}{\psi_q},{\bf z_{q+1}}\right)\widetilde{f}\left(\frac{x-\sum_{j\geq q+1}\psi_j z_j}{\psi_q}\vert {\bf z_{q+1}}\right)\\&&y dx dy d\phi dz_3\cdots dz_{q-1} d\nu({\bf z_{q+1}})\\
\end{eqnarray*}
Then we get 
\begin{eqnarray*}
\omega(x,y)&&=y\int \zeta_1\left(\sqrt{y^2-k_3({\bf z_3})}\cos\phi, \sqrt{y^2-k_3({\bf z_3})}\cos\phi,z_3,\ldots,z_{q-1},\frac{x-\sum_{j\geq q+1}\psi_j z_j}{\psi_q},{\bf z_{q+1}}\right)\\
&&\zeta_2\left(z_3,\ldots,z_{q-1},\frac{x-\sum_{j\geq q+1}\psi_j z_j}{\psi_q},{\bf z_{q+1}}\right)\\
&&\widetilde{f}\left(\frac{x-\sum_{j\geq q+1}\psi_j z_j}{\psi_q}\vert {\bf z_{q+1}}\right)d\phi dz_3\cdots dz_{q-1} d\nu({\bf z_{q+1}}).          \end{eqnarray*}
We deduce that there exists $C>0$ such that $\omega(x,y)\leq C y$. $\square$   
\end{enumerate}

\subsection{Auxiliary Lemmas}
This subsection presents two auxiliary Lemma which assert that under our assumptions, truncated versions $\overline{m}_t$ and $\overline{\sigma}_t$ of $m_t$ and $\sigma_t$ have no effect in the asymptotic expansion of our estimator. Moreover, $m_t$ and $\sigma_t$ can be replaced with their $\ell-$dependent approximations, provided $\ell=\ell_n$ grows at an arbitrary small power of $n$.

We first observe that if the kernel $K$ is Lipschitz continuous and bounded, we have for any $q\in (0,1)$, 
$$\left|K(x)-K(y)\right|\leq \left(2\Vert K\Vert_{\infty}\right)^{1-q}\mbox{Lip}(K)^q\vert x-y\vert^q.$$
In particular, $K$ is H\"{o}lder continuous with exponent $q$. This fact will be used in the two following
lemmas.  
In what follows, we set for $1\leq i,j\leq n$ and $(v,\theta)\in I\times \Theta_{0,\epsilon}$,
$$\mathcal{L}_{v,ij}(\theta)=L_{v,i}(\theta)-\ep_j(\theta),\quad \mathcal{L}_{v,ij\ell}(\theta)=L_{v,i\ell}(\theta)-\ep_{j\ell}(\theta),\quad \overline{\mathcal{L}}_{v,ij}(\theta)=\overline{L}_{v,i}(\theta)-\overline{\ep}_j(\theta).$$
Here $L_{v,ij\ell}(\theta)=\frac{v-m_{i\ell}(\theta)}{\sigma_{i\ell}(\theta)}$ and $\ep_{j\ell}(\theta)=\frac{X_{j\ell}-m_{j\ell}(\theta)}{\sigma_{j\ell}(\theta)}$. We also set
$$\widetilde{f}_X(v)=\frac{1}{n^2}\sum_{i,j=1}^n\frac{1}{\sigma_i(\hat{\theta})}K_b\left[\mathcal{L}_{v,ij}(\hat{\theta})\right].$$
\begin{lem}\label{approx1}
Assume that there exists $\delta\in (0,1)$ such that  $nb^{2+\delta}\rightarrow\infty$. Then we have $\sup_{v\in I}\left|\hat{f}_X(v)-\widetilde{f}_X(v)\right|=o_{\P}\left(\frac{1}{\sqrt{n}}\right)$. 
\end{lem}

\paragraph{Proof of Lemma \ref{approx1}}
Since $\sqrt{n}\left(\hat{\theta}-\theta_0\right)=O_{\P}(1)$, it is enough to prove that
$$\frac{1}{n^2}\sum_{i,j=1}^n\sup_{v\in I,\theta\in \Theta_{0,\epsilon}}\left|\frac{1}{\sigma_i(\theta)}K_b\left[\mathcal{L}_{v,ij}(\theta)\right]-\frac{1}{\sigma_{i\ell}(\theta)}K_b\left[\overline{\mathcal{L}}_{v,ij}(\theta)\right]\right|=o_{\P}\left(n^{-1/2}\right),$$
where $\epsilon>0$ is defined in assumption {\bf A2}.
Let $q$ be a positive real number such that $2q\leq \min(\delta,s)$. Using assumption {\bf A2}, we have 
$$\max_{\theta\in\Theta}\left|\frac{1}{\sigma_i(\theta)}-\frac{1}{\overline{\sigma}_i(\theta)}\right|\leq C \max_{\theta\in \Theta}\left|\sigma_i^2(\theta)-\overline{\sigma}_i^2(\theta)\right|^q.$$
One can choose for instance $C=2^{-q}\gamma^{-1-2q}$.
Moreover, setting $q=s/2$, we have 
\begin{eqnarray*}
&&\left|K_b\left[\mathcal{L}_{v,ij}(\theta)\right]-K_b\left[\overline{\mathcal{L}}_{v,ij}(\theta)\right]\right|\\
&\leq & Cb^{-1-q}\left|\mathcal{L}_{v,ij}(\theta)-\overline{\mathcal{L}}_{v,ij}(\theta)\right|^q\\
&\leq &C b^{-1-q}\left[\left|m_i(\theta)-\overline{m}_i(\theta)\right|^q+\left|m_i(\theta)\right|^q\cdot\left|\sigma^2_i(\theta)-\overline{\sigma}^2_i(\theta)\right|^q\right]\\
&+&Cb^{-1-q}\left[\left|m_j(\theta)-\overline{m}_j(\theta)\right|^q+\left|\sigma^2_j(\theta)-\overline{\sigma}^2_j(\theta)\right|^q\cdot\left|X_j-m_j(\theta)\right|^q\right].
\end{eqnarray*}
Using assumption {\bf A2} and the Cauchy-Schwarz inequality, we deduce that 
$$\sup_{v\in I}\left\|\hat{f}(v)-\widetilde{f}(v)\right\|=O_{\P}\left(\frac{1}{nb^{1+q}}\sum_{i=1}^na^{\frac{iq}{s}}\right)=O_{\P}\left(\frac{1}{nb^{1+q}}\right).\square$$

\begin{lem}\label{approx2}
Assume that $\ell=n^t$ with $t>0$ and that $nb\rightarrow \infty$. Then if
$$\widetilde{f}_{\ell}(v)=\frac{1}{n^2}\sum_{i,j=1}^n\frac{1}{\sigma_{i\ell}(\hat{\theta})}K_b\left[L_{v,i\ell}(\hat{\theta})-\ep_{j\ell}(\hat{\theta})\right].$$
We have $\sup_{v\in I}\left|\widetilde{f}_X(v)-\widetilde{f}_{\ell}(v)\right|=o_{\P}\left(n^{-1/2}\right)$.
\end{lem}

\paragraph{Proof of Lemma \ref{approx2}}
Since $\sqrt{n}\left(\hat{\theta}-\theta_0\right)=O_{\P}(1)$, it is enough to prove that
$$\frac{1}{n^2}\sum_{i,j=1}^n\sup_{v\in I,\theta\in \Theta_{0,\epsilon}}\left|\frac{1}{\sigma_i(\theta)}K_b\left[\mathcal{L}_{v,ij}(\theta)\right]-\frac{1}{\sigma_{i\ell}(\theta)}K_b\left[\mathcal{L}_{v,ij\ell}(\theta)\right]\right|=o_{\P}\left(n^{-1/2}\right),$$
where $\epsilon>0$ is defined in assumption {\bf A2}.
\begin{itemize}
\item
As in the proof of Lemma \ref{approx1}, we use assumption {\bf A2} to get 
$$\left|\frac{1}{\sigma_i(\theta)}-\frac{1}{\sigma_{i\ell}(\theta)}\right|\leq C \max_{\theta\in \Theta}\left|\sigma_i^2(\theta)-\sigma_{i\ell}^2(\theta)\right|^s.$$
\item
Using the previous point, we have also for $q=s/2$,
\begin{eqnarray*}
&&\sup_{v\in I}\left|K_b\left[\mathcal{L}_{v,ij}(\theta)\right]-K_b\left[\mathcal{L}_{v,ij\ell}(\theta)\right]\right|\\
&\leq & Cb^{-1-q}\sup_{v\in I}\left|\mathcal{L}_{v,ij}(\theta)-\mathcal{L}_{v,ij\ell}(\theta)\right|^q\\
&\leq &C b^{-1-q}\left[\left|m_i(\theta)-m_{i\ell}(\theta)\right|^q+\left|m_i(\theta)\right|^q\cdot\left|\sigma^2_i(\theta)-\sigma^2_{i\ell}(\theta)\right|^q\right]\\
&+&Cb^{-1-q}\left[\left|m_j(\theta)-m_{j\ell}(\theta)\right|^q+\left|X_j-X_{j\ell}\right|^q
+\left|\sigma^2_j(\theta)-\sigma^2_{j\ell}(\theta)\right|^q\cdot\left|X_j-m_j(\theta)\right|^q\right].
\end{eqnarray*}
Using assumption {\bf A2}, we get 
$$\E\sup_{\theta\in\Theta_{0,\epsilon}, v\in I}\left|K_b\left[\mathcal{L}_{v,ij}(\theta)\right]-K_b\left[\mathcal{L}_{v,ij\ell}(\theta)\right]\right|\leq C \frac{a^{\ell/2}}{b^{1+q}}.$$
\end{itemize}
Using the two previous points and assumption {\bf A2}, we get 
$$\frac{1}{n^2}\sum_{i,j=1}^n\sup_{v\in I,\theta\in \Theta_{0,\epsilon}}\left|\frac{1}{\sigma_i(\theta)}K_b\left[\mathcal{L}_{v,ij}(\theta)\right]-\frac{1}{\sigma_{i\ell}(\theta)}K_b\left[\mathcal{L}_{v,ij\ell}(\theta)\right]\right|=O_{\P}\left(\frac{a^{\ell/2}}{b^{1+q}}\right).$$
Then the result follows from the conditions $\ell\sim n^t$, $0<a<1$ and $nb\rightarrow \infty$.$\square$

\subsection{Residual process regularity and an auxiliary lemma}
In this subsection, we provide auxiliary lemmas which give regularity properties of the density of the residual process as well as 
some moment conditions for its $\ell-$dependent approximation. Then we will state Lemma \ref{propres2} which will be  
required for the proof of Theorem $2$. For $\theta\in \Theta_{0,\epsilon}$, the density of $\ep_t(\theta)$ will be denoted by $f_{\theta}$.

We have the expression 
\begin{equation}\label{exprime}
f_{\theta}(w)=\E\left[\frac{\sigma_j(\theta)}{\sigma_j(\theta_0)}f_{\ep}\left(\frac{\sigma_j(\theta)}{\sigma_j(\theta_0)}w+\frac{m_j(\theta)-m_j(\theta_0)}{\sigma_j(\theta_0)}\right)\right].
\end{equation}  

The following lemma is given without proof. The assertions given below are mainly a consequence of the Lebesgue theorem for the derivative of an integral depending on a parameter. 
\begin{lem}\label{propres}
Assume that assumptions {\bf A3} and {\bf A5} hold true. 
\begin{enumerate}
\item
We have $\sup_{\theta\in\Theta_{0,\epsilon},w\in \R}f_{\theta}(w)<\infty$. For $\theta\in \Theta_{0,\epsilon}$, the function $w\mapsto f_{\theta}(w)$ has a derivative $f'_{\theta}$ such that $\sup_{\theta\in \Theta_{0,\epsilon},w\in \R}\left|f'_{\theta}(w)\right|<\infty$. Moreover, there exists a number $C>0$ not depending on $\theta,w$ such that
$$\left|f'_{\theta}(w)-f'_{\theta_0}(w)\right|\leq C\left(1+|w|\right)\Vert \theta-\theta_0\Vert.$$ 
\item
For each $w\in\R$, the function $\theta\mapsto f_{\theta}(w)$ is two times differentiable. 
Moreover, there exists a number $C>0$ not depending on $w,\theta$ such that 
$$\left|\dot{f}_{\theta}(w)-\dot{f}_{\theta_0}(w)\right|\leq C\left(1+w^2\right)\Vert \theta-\theta_0\Vert,\quad \left|\dot{f}_{\theta}(w)\right|\leq C\left(1+\vert w\vert\right)$$
and 
$$\left|\dot{f}_{\theta_0}(w_1)-\dot{f}_{\theta_0}(w_2)\right|\leq C\left(1+\vert w_1\vert\right)\cdot\vert w_1-w_2\vert.$$
\item
There exist two constants $C_1$ and $C_2$ such that
$$\left\vert f_{\theta}(w)-f_{\theta_0}(w)\right\vert\leq \left(C_1+C_2|w|\right)\Vert\theta-\theta_0\Vert.$$
\item
If $F:\R\rightarrow \R$ a function continuously differentiable and with a compact support. Then, for $\theta\in B(\theta_0,\epsilon)$, we have 
$$\E\left[F\left(\ep_t(\theta)\right)\right]=\int F(w)f_{\theta}(w)dw,\quad \E\left[\dot{\ep}_t(\theta)F'\left(\ep_t(\theta)\right)\right]=
\int F(w) \dot{f}_{\theta}(w)dw.$$
\end{enumerate}
\end{lem}

The next lemma is given without proof because it results from simple computations.
\begin{lem}\label{propres3}
Assume that assumption {\bf A3} holds true. We set $U_j(\theta)=\frac{m_j(\theta_0)-m_j(\theta)}{\sigma_j(\theta)}$ and $V_j(\theta)=\frac{\sigma_j(\theta_0)}{\sigma_j(\theta)}$. Note that $\ep_j(\theta)=U_j(\theta)+V_j(\theta)\ep_j$.
\begin{enumerate}
\item
We have $\E\nor{\dot{U_j}}^2<\infty$, $\E\nor{\frac{U_j\dot{V_j}}{V_j}}^2<\infty$, $\E\nor{\frac{\dot{V_j}}{V_j}}^2<\infty$, $\E\nor{\dot{\ep}_j}^2<\infty$ and $\E\nor{\ddot{\ep}_{j}}<\infty$.
\item
We have  $\E\left[\sup_{v\in I}\nor{L_{v,i}}^2\right]<\infty$, $\E\left[\sup_{v\in I}\nor{\dot{L}_{v,i}}^2\right]<\infty$ and $\E\left[\sup_{v\in I}\nor{\ddot{L}_{v,i}}\right]<\infty$.
\item
Let $\eta$ be a positive number. There exists a positive real number $C$ such that
$$\sup_{\substack{\vert v-w\vert\leq \eta\\ \Vert\theta-\zeta\Vert\leq \eta}}\left|
L_{v,i}(\theta)-L_{w,i}(\zeta)\right|\leq C\eta\left(1+sup_{v\in I}\nor{\dot{L}_{v,i}}\right)$$
and 
$$\sup_{\substack{\vert v-w\vert\leq \eta\\ \Vert\theta-\zeta\Vert\leq \eta}}\left\|
\dot{L}_{v,i}(\theta)-\dot{L}_{w,i}(\zeta)\right\|\leq C\eta\left(\nor{\dot{\sigma}_{i},\sigma_{i}}+sup_{v\in I}\nor{\ddot{L}_{v,i}}\right).$$

\end{enumerate}
\end{lem}

Now we state a lemma which will be useful for studying the uniform convergence of some sums of martingale differences.
We set $n'=k\ell$ where $k$ is the integer part of $n/\ell$. Moreover, we set 
$$\mathcal{I}_n=\left\{1\leq i,j\leq n': i\leq j-\ell\mbox{ or } i\geq i+\ell\right\}.$$
\begin{lem}\label{propres2}
Assume that assumptions {\bf A3}and {\bf A5} hold true. 
\begin{enumerate}
\item
We set $B_{v,ij}(\theta)=\dot{A}_{v,ij}(\theta)-\E_{Y_{j\ell}}\left[\dot{A}_{v,ij}(\theta)\right]$. 
Then we have 
$$\sum_{(i,j)\in \mathcal{I}_n}\sup_{\substack{v\in I\\\theta\in B(\theta_0,\epsilon)}}\E_{Y_{j\ell}}\Vert B_{v,ij}(\theta)\Vert^2=O_{\P}\left(\frac{n^2}{b^3}\right).$$
\item
We have 
$$\max_{1\leq j\leq n}\sum_{\substack{i: (i,j)\in\mathcal{I}_n\\1\leq i\leq n}}^n\sup_{v\in I}\left|B_{v,ij}\right|_{\infty,\epsilon}=O_{\P}\left(\frac{n^{4/3}}{b^2}\right).$$  
\item
Let $\eta$ be a positive number. We have
$$\sum_{(i,j)\in\mathcal{I}_n}\sup_{\substack{\vert v-w\vert\leq \eta\\\Vert\theta-\zeta\Vert\leq\eta}}\left\Vert B_{v,ij}(\theta)-B_{w,ij}(\zeta)\right\Vert=O_{\P}\left(\frac{\eta n^2}{b^3}\right).$$
\end{enumerate}
\end{lem}
\paragraph{Proof of Lemma \ref{propres2}.}
We only prove the result for the conditionally heteroscedastic case, the homoscedastic uses similar arguments and is simpler.
\begin{enumerate}
\item
It is only necessary to prove that 
$$\sum_{(i,j)\in \mathcal{I}_n}\sup_{\substack{v\in I\\\theta\in B(\theta_0,\epsilon)}}\E_{Y_{j\ell}}\Vert \dot{A}_{v,ij}(\theta)\Vert^2=O_{\P}\left(\frac{n^2}{b^3}\right).$$
We recall that 
$$\dot{A}_{v,ij}(\theta)=\dot{\sigma_{i\ell}^{-1}}(\theta)K_b\left[L_{v,i\ell}(\theta)-\ep_{j\ell}(\theta)\right]+\sigma_{i\ell}^{-1}(\theta)\left[\dot{L}_{v,i\ell}(\theta)-\dot{\ep}_{j\ell}(\theta)\right]K_b'\left[L_{v,i\ell}(\theta)-\ep_{j\ell}(\theta)\right].$$
We define 
$$U_{j\ell}(\theta)=\frac{m_{j\ell}(\theta_0)-m_{j\ell}(\theta)}{\sigma_{j\ell}(\theta)},\quad V_{j\ell}(\theta)=\frac{\sigma_{j\ell}(\theta_0)}{\sigma_{j\ell}(\theta)}.$$
Then $\ep_{j\ell}(\theta)=U_{j\ell}(\theta)+V_{j\ell}(\theta)\ep_j$. Moreover 
\begin{eqnarray*}
&&\E_{\ep}\left[\Vert\dot{\ep}_{j\ell}(\theta)\Vert^2\cdot \left\vert K'_b\left(L_{v,i\ell}(\theta)-\ep_{j\ell}(\theta)\right)\right\vert^2\right]\\
&=& \int \frac{1}{b^3}\Vert \dot{U}_{j\ell}(\theta)+\frac{\dot{V}_{j\ell}(\theta)}{V_{j\ell}(\theta)}\left[L_{v,i\ell}(\theta)-U_{j\ell}(\theta)-bw\right]\Vert^2\cdot
\left\vert K'(w)\right\vert^2 f_{\ep}\left(\frac{L_{v,i\ell}(\theta)-U_{j\ell}(\theta)-bw}{V_{j\ell}(\theta)}\right)dw\\
&\leq& \frac{C}{b^3}\left[1+\nor{\dot{U}_{j\ell}}^2+\nor{\frac{\dot{V}_{j\ell}}{V_{j\ell}}}^2+\sup_{v\in I}\nor{L_{v,i\ell}}^2\cdot\nor{\frac{U_{j\ell}\dot{V}_{j\ell}}{V_{j\ell}}}^2\right].
\end{eqnarray*}
Then we get 
\begin{eqnarray*}
&&\E_{Y_{j\ell}}\left[\Vert \dot{A}_{v,ij}(\theta)\Vert^2\right]\\
&\leq & \frac{C}{b^3}\left[1+\nor{\dot{\sigma^2_{i\ell}},\sigma^2_{i\ell}}^2+\sup_{v\in I} \nor{\dot{L}_{v,i\ell}}^2+\E\nor{\dot{U}_{j\ell}}^2+\E\nor{\frac{\dot{V}_{j\ell}}{V_{j\ell}}}^2+\sup_{v\in I}\nor{L_{v,i\ell}}^2\cdot\E\nor{\frac{U_{j\ell}\dot{V}_{j\ell}}{V_{j\ell}}}^2\right].
\end{eqnarray*}
The result follows from assumption {\bf A3} and Lemma \ref{propres3}.
\item
Since 
$$\dot{\ep}_{j\ell}(\theta)=-\frac{\dot{m}_{j\ell}(\theta)}{\sigma_{j\ell}(\theta)}-\frac{\dot{\sigma}_{j\ell}(\theta)}{\sigma_{j\ell}(\theta)}\ep_{j\ell}(\theta),$$
we have, using the compact support of the kernel $K$, and the equality $\frac{\dot{\sigma}_j}{\sigma_j}=\frac{\dot{\sigma^2_j}}{2\sigma^2_j}$, 
$$\left\Vert \dot{\ep}_{j\ell}(\theta)\right\Vert\cdot\left\vert K_b'\left(L_{v,i\ell}(\theta)-\ep_{j\ell}(\theta)\right)\right\vert\leq \frac{C}{b^2}\left[\nor{\dot{m}_{j\ell},\sigma_{j\ell}}+\nor{\dot{\sigma^2_{j\ell}},\sigma^2_{j\ell}}\cdot\left[1+\sup_{v\in I}\nor{L_{v,i\ell}}\right]\right].$$
Then we conclude that  
\begin{eqnarray*}
&&\sup_{v\in I}\left\Vert \dot{A}_{v,ij}\right\|\\
&\leq& \frac{C}{b^2}\left[\nor{\dot{\sigma^2_{i\ell}},\sigma^2_{i\ell}}+\sup_{v\in I}\nor{\dot{L}_{v,i\ell}}+\nor{\dot{m}_{j\ell},\sigma_{j\ell}}+
\nor{\dot{\sigma^2_{j\ell}},\sigma^2_{j\ell}}\cdot\left[1+\sup_{v\in I}\nor{L_{v,i\ell}}\right]\right].
\end{eqnarray*}
From the assumption {\bf A3}, we have $\max_{1\leq j\leq n}\nor{\dot{\sigma^2_{j\ell}},\sigma^2_{j\ell}}=O_{\P}\left(n^{1/3}\right)$
and $\max_{1\leq j\leq n}\nor{\dot{m}_{j\ell},\sigma_{j\ell}}=O_{\P}\left(n^{1/3}\right)$. Then the result follows from the point $2$ of Lemma \ref{propres3}. 
\item
If $(v,w)\in I^2$ and $(\theta,\zeta)\in \Theta_{0,\epsilon}$ are such that
$\vert v-w\vert\leq\eta$ and $\Vert\theta-\zeta\Vert\leq \eta$, some basic computations lead to the inequality 
\begin{eqnarray*}
\left\vert \dot{A}_{v,ij}(\theta)-\dot{A}_{w,ij}(\zeta)\right\vert&\leq& \frac{C\eta}{b^3}\left[1+\nor{\ddot{\sigma^2_{i\ell}},\sigma^2_{i\ell}}
+\nor{\dot{\sigma^2_{i\ell}},\sigma^2_{i\ell}}^2\right]\\
&+& \frac{C\eta}{b^3}\left[\sup_{v\in I}\nor{\dot{L}_{v,i\ell}}^2+\nor{\dot{\ep}_{j\ell}}^2+\sup_{v\in I}\nor{\ddot{L}_{v,i\ell}}+\nor{\ddot{\ep}_{j\ell}}\right].
\end{eqnarray*}
Then the result follows from assumption {\bf A3} and Lemma \ref{propres3}.$\square$
\end{enumerate}

\bibliographystyle{plainnat}
\bibliography{bibDensity}

\end{document}